\documentclass[leqno,twoside,11pt]{amsart}

\usepackage{latexsym,esint,url}
\usepackage{color}

\theoremstyle{plain}
\def\endproof{\hspace*{\fill}\mbox{\ \rule{.1in}{.1in}}\medskip }

\newtheorem{theorem}{Theorem}[section]
\newtheorem{corollary}[theorem]{Corollary}
\newtheorem{lemma}[theorem]{Lemma}

\theoremstyle{definition}

\newtheorem{remark}[theorem]{Remark}

\numberwithin{equation}{section}
\numberwithin{figure}{section}

\newcommand{\e}{\varepsilon}

\begin{document}

\title[The matching property of infinitesimal isometries for elliptic shells]
{The matching property of infinitesimal isometries \\ on elliptic surfaces\\ 
and elasticity of thin shells}
\author{Marta Lewicka, Maria Giovanna Mora and Mohammad Reza Pakzad}
\address{Marta Lewicka, University of Minnesota, Department of Mathematics, 
206 Church St. S.E., Minneapolis, MN 55455}
\address{Maria Giovanna Mora, Scuola Internazionale Superiore di Studi Avanzati,
via Beirut 2-4, 34014 Trieste, Italy}
\address{Mohammad Reza Pakzad, University of Pittsburgh, Department of Mathematics, 
139 University Place, Pittsburgh, PA 15260}
\email{lewicka@math.umn.edu, mora@sissa.it, pakzad@pitt.edu}

\date{\today}
\begin{abstract} 
Using the notion of $\Gamma$-convergence, we discuss the limiting behavior of the 3d nonlinear 
elastic energy for thin elliptic shells, as their thickness $h$ converges to
zero, under the assumption that the elastic energy of deformations scales 
like $h^\beta$ with $2<\beta<4$. We establish that, for the given scaling 
regime, the limiting theory reduces to the linear pure bending. 
Two major ingredients of the proofs are: the density of smooth 
infinitesimal isometries in the space of $W^{2,2}$ first order infinitesimal
isometries, and a result on matching smooth infinitesimal isometries with 
exact isometric immersions on smooth elliptic surfaces. 
\end{abstract}

\maketitle
\tableofcontents

\section{Introduction} 

In this paper we continue the rigorous derivation of shell theories by three-dimensional
nonlinear elasticity. We follow our earlier work \cite{lemopa1}
in which we derived the generalization of the von K\'arm\'an theory (introduced in the setting
of plates and justified in \cite{FJMhier} through $\Gamma$-convergence) to shells
with mid-surfaces of arbitrary geometry, and \cite{lemopa2} where we treated shells with
variable thickness.  From the mathematical elasticity point of view, the present paper
generalizes the limiting theory which in \cite{FJMhier} corresponds to the linearized
Kirchhoff model for plates. See also \cite{FJMgeo, FJMM_cr, LR1, LeD-Rao, CM05} for rigorous derivations
of other theories from the same point of view. We refer to \cite{ciarbookvol3} for a classical
discussion of such theories, using the asymptotic expansion, although we point out that none
of the scaling regimes in \cite{lemopa1} or in this paper have been discussed in the general
context of shells before.

\medskip

Consider a $2$-dimensional surface $S$ embedded in $\mathbb{R}^3$, which is
compact, connected, oriented, of class $\mathcal{C}^{1,1}$,
and with boundary $\partial S$ being the union of finitely many 
(possibly none) Lipschitz curves.
A family $\{S^h\}_{h>0}$ of shells of small thickness $h$
around $S$ is given through:
$$S^h = \{z=x + t\vec n(x); ~ x\in S, ~ -h/2< t < h/2\},\qquad 0<h<h_0.$$
By $\vec n(x)$ we denote the unit normal to $S$, by $T_x S$ the tangent space,
and $\Pi(x) = \nabla \vec n(x)$ is the shape operator on $S$
(or equivalently: its negative second fundamental form). 
The projection onto $S$ along $\vec n$ will be denoted by $\pi$.
We will assume that $h<h_0$, with $h_0$ sufficiently small
to have $\pi$ well defined on each $S^h$.

To a deformation $u\in W^{1,2}(S,\mathbb{R}^3)$ we associate its elastic
energy (scaled per unit thickness):
\begin{equation}\label{elastic-En}
E^h(u) = \frac{1}{h}\int_{S^h} W(\nabla u).
\end{equation} 
Here, the stored energy density $W:\mathbb{R}^{3\times 3}\longrightarrow [0,\infty]$ 
is assumed to be $\mathcal{C}^2$ in a neighborhood of $SO(3)$, and to satisfy
the following normalization, frame indifference and nondegeneracy conditions:
\begin{equation*}\label{assump-intro}
\begin{split}
\forall F\in\mathbb{R}^{3\times 3}\quad \forall R\in SO(3)\qquad
& W(R) = 0, \quad W(RF) = W(F), \\
& W(F)\geq C \mathrm{dist}^2(F, SO(3)) 
\end{split}
\end{equation*}
(with a uniform constant $C>0$). In the study of the elastic properties of
thin shells $S^h$, a crucial step is to describe the limiting behavior, as
$h\to 0$, of minimizers $u^h$ to the total energy functional:
\begin{equation}\label{total-intro}
J^h(u) = E^h(u) - \frac{1}{h}\int_{S^h} f^hu,
\end{equation} 
subject to applied forces $f^h$. It can be shown that if the forces $f^h$ scale like
$h^\alpha$, then $E^h(u^h)\sim h^\beta$ where
$\beta= \alpha$ if $0 \le \alpha \le 2$ and $\beta
= 2\alpha -2$ if $\alpha > 2$. The main part of the analysis consists 
therefore of characterizing the
limiting behavior of the scaled energy functionals $1/h^\beta E^h$, or more generally, that of
$1/e^h E^h$, where $e^h$ is a given sequence of positive numbers obeying a prescribed scaling
law.  Throughout this paper we shall assume that:
\begin{equation}\label{scaling-intro}
0<\lim_{h\to 0} e^h/h^\beta < +\infty, \qquad \mbox{ for some }~~ 2<\beta<4.
\end{equation}

The first result in this framework is due to LeDret and Raoult \cite{LeD-Rao}, who studied
the scaling $\beta=0$. It leads to a membrane shell model, with energy depending only 
on stretching and shearing of the mid-surface. The case $\beta=2$ has been
analyzed in \cite{FJMM_cr} and it corresponds to the geometrically nonlinear
bending theory, where the only admissible deformations are 
the isometries of the mid-surface, while the energy expresses the total change
of curvature produced by the deformation.

In the recent paper \cite{lemopa1} the limiting model has been identified for
the range of scalings $\beta\geq 4$. In these cases, the admissible
deformations $u$ are only those which are close to a rigid motion
$R$ and whose first order term in the expansion of $u-R$ with respect to $h$
is given by $RV$, where $V$ is an element of the 
class $\mathcal V$ of {\em infinitesimal isometries} on $S$ \cite{Spivak}. 
More precisely, $\mathcal{V}$ consists of vector fields
$V\in W^{2,2}(S,\mathbb{R}^3)$ for whom there 
exists  a matrix field $A\in W^{1,2}(S,\mathbb{R}^{3\times 3})$ so that:
\begin{equation}\label{Adef-intro} 
\partial_\tau V(x) = A(x)\tau \quad \mbox{and} \quad  A(x)^T= -A(x) \qquad 
\forall {\rm{a.e.}} \,\, x\in S \quad \forall \tau\in T_x S.
\end{equation}
The class $\mathcal V$ naturally plays crucial role in the analysis of shells,
and its members are identified by the following geometric property: $V$ is a
(first order) infinitesimal isometry if the change of metric induced by the
deformation $\mbox{id} + \e V$ is at most of order $\e^2$.

For $\beta>4$ the limiting energy is given only by a bending term, that is 
the first order change in the second fundamental 
form of $S$, produced by $V$:
\begin{equation}\label{I-intro}
I(V)= \frac{1}{24} \int_S 
\mathcal{Q}_2\left(x,(\nabla(A\vec{n}) - A\Pi)_{tan}\right)~\mbox{d}x,
\qquad \forall V\in\mathcal{V},
\end{equation} 
and corresponds to the linear pure bending theory
derived in \cite{ciarbookvol3} from linearized elasticity.  

In (\ref{I-intro}), the quadratic 
forms $\mathcal{Q}_2(x,\cdot)$ are defined as follows:
$$ \mathcal{Q}_2(x, F_{tan}) = \min\{\mathcal{Q}_3(\tilde F); ~~ (\tilde F - F)_{tan} = 0\}, 
\qquad \mathcal{Q}_3(F) = D^2 W(\mbox{Id})(F,F).$$
The form $\mathcal{Q}_3$ is defined for all $F\in\mathbb{R}^{3\times 3}$, 
while $\mathcal{Q}_2(x,\cdot)$,  for a given $x\in S$ is defined on tangential 
minors $F_{tan}$ of such matrices.   Recall that the tangent 
space to $SO(3)$ at $\mbox{Id}$ is $so(3)$. As a consequence,  
both forms depend only on the symmetric parts of their arguments
and are positive definite on the space of symmetric matrices \cite{FJMgeo}. 

For $\beta=4$ the $\Gamma$-limit which turns out to be the generalization of 
the von K\'arm\'an functional \cite{FJMhier} to shells, 
contains also a stretching term, measuring the second order change in the
metric of $S$:
\begin{equation*}\label{vonKarman}
\tilde I(V,B_{tan})= \frac{1}{2} 
\int_S \mathcal{Q}_2\left(x,B_{tan} - \frac{1}{2} (A^2)_{tan}\right)
+ \frac{1}{24} \int_S \mathcal{Q}_2\left(x,(\nabla(A\vec n) - A\Pi)_{tan}\right).
\end{equation*}
It involves a symmetric matrix field $B_{tan}$ belonging to the {\em finite strain space}: 
$$ \mathcal{B} = \Big\{L^2 - \lim_{h\to 0}\mathrm{sym }\nabla w^h; 
~~ w^h\in W^{1,2}(S,\mathbb{R}^3)\Big\}. $$ 
The space $\mathcal{B}$ emerges as well in the context of  linear elasticity
and ill-inhibited surfaces \cite{sanchez, GSP}.

It was further shown in \cite{lemopa1} that for a certain class of surfaces, 
referred to as {\it approximately robust surfaces}, the limiting 
theory for $\beta=4$ reduces to the purely linear bending functional (\ref{I-intro}). 
Strictly elliptic (or convex) $\mathcal{C}^{2,1}$ surfaces  
happen to belong to this class \cite{lemopa1}.

\medskip

Here we focus on the range of scalings $2<\beta<4$. 
Roughly speaking, we look for an intermediate theory between the limit
theories corresponding to $\beta=2$ and $\beta\ge 4$. 
On one hand, modulo a rigid motion, the deformation of the mid-surface must 
look like $\mbox{id} + \e V$ up to its first order of expansion. 
On the other hand, the closer $\beta$ is to $2$, 
the closer the mid-surface deformation must be to an exact isometry of $S$. 
To overcome this apparent disparity between first order infinitesimal
isometries and exact isometries in this context, one must study the conditions 
under which, given $V \in {\mathcal V}$, one can construct an exact isometry 
of the form $\mbox{id} + \e V + \e^2 w_\e$, with equibounded $w_\e$. This is what 
we refer to as the {\em matching property} of an infinitesimal isometry. 
A key  question is hence whether the infinitesimal isometry $V$ obtained as
above from the limit deformation $u$ enjoys the matching property for 
the energy scaling range $2<\beta<4$.

If $S\subset \mathbb R^2$ represents a plate, the above questions have been
answered  in \cite{FJMhier}. In this case:
\begin{itemize}
\item[(1)] The limit displacement $V$ must necessarily belong to the space of
  second order infinitesimal isometries:
$$ \mathcal V_2:= \{ V\in {\mathcal V};~~  (A^2)_{tan} \in {\mathcal B} \},$$ 
where the matrix field $A$ is as in (\ref{Adef-intro}).  
\item[(2)] Any Lipschitz second order isometry $V\in\mathcal{V}_2$
satisfies the matching property. 
\end{itemize}  
Combining these two facts with the density of Lipschitz second order 
infinitesimal isometries in ${\mathcal V_2}$ for a plate \cite{MuPa2},  
one concludes through careful $\Gamma$-convergence arguments that the 
limiting plate theory is given by the functional (\ref{I-intro}) over the
nonlinear space ${\mathcal V_2}$.  
Note that, for a plate, $V\in {\mathcal V_2}$ means that there exists an 
in-plane displacement $w\in W^{1,2} (S, \mathbb R^2)$ 
such that the change of metric due to $\mbox{id} + \e V + \e^2 w$ is of order $\e^3$.
Also, in this case, an equivalent analytic characterization for 
$V=(V^1, V^2, V^3) \in {\mathcal V_2}$ is given by 
$(V^1, V^2) =  (-\omega y, \omega x) + (b_1, b_2)$ and $\det \nabla^2 V^3 =0$. 

\medskip

Towards analyzing more general surfaces $S$, we first derive the matching
property and the corresponding density of isometries, for elliptic surfaces.
We say that $S$ is elliptic if its shape operator $\Pi$ is strictly positive 
(or strictly negative) definite up to the
boundary. Hence, without loss of generality we have:
\begin{equation}\label{elliptic-intro} 
\forall x\in \bar S \quad \forall \tau\in T_xS \qquad 
\frac{1}{C}|\tau|^2 \leq \Big(\Pi(x)\tau\Big)\cdot\tau \leq C|\tau|^2.
\end{equation}  
The main result of this paper states that for elliptic surfaces of
sufficient regularity, the $\Gamma$-limit of 
the nonlinear elastic energy (\ref{elastic-En}) for the scaling regime
$2<\beta<4$ is still given by the energy functional (\ref{I-intro}) 
over the linear space ${\mathcal V}$.  
The main reason is the fact that for an elliptic surface, all sufficiently smooth infinitesimal 
isometries satisfy the matching property. This implies a qualitatively
different theory from the nonlinear theory obtained for plates in this regime.
Our main theorem about the matching property of infinitesimal isometries is the following:
\begin{theorem}\label{th_exact-intro}
Let $S$ be elliptic (that is condition \eqref{elliptic-intro} is satisfied). 
Moreover assume that  $S$ is homeomorphic to a disk and that for some $\alpha
>0$,  $S$ and $\partial S$ are of class $\mathcal{C}^{3,\alpha}$. 
Given $V\in\mathcal{V}\cap\mathcal{C}^{2,\alpha}(\bar S)$, there exists 
a sequence $w_h:\bar S\longrightarrow \mathbb{R}^3$, equibounded in 
$\mathcal{C}^{2,\alpha}(\bar S)$, and such that for all small $h>0$ the map
$u_h = \mathrm{id} +hV + h^2w_h$ is an (exact) isometry. 
\end{theorem}  

We apply this result in section \ref{sec_recovery},
to construct the recovery sequence necessary for establishing the upper bound in the context 
of the $\Gamma$-convergence (see Theorem \ref{thlimsup-intro}). 
Clearly, Theorem \ref{th_exact-intro} is not sufficient for this purpose as
the elements of $\mathcal V$ are only $W^{2,2}$ regular. Indeed, in 
most $\Gamma$-convergence results, a key step is to prove density of suitable 
more regular mappings in the space of admissible mappings for the limit problem. Results 
in this direction, for Sobolev spaces of isometries and infinitesimal
isometries of flat regions, have already been 
proved and applied  in the context of derivation of plate theories. 
The interested reader can refer to \cite{Pak, MuPa2, Ho1, Ho2} 
for statements of these density theorems and their applications in \cite{FJMhier, CD}.

In general, even though $\mathcal{V}$ is
a linear space, and assuming $S$ to be  $\mathcal{C}^\infty$, the usual mollification 
techniques do not guarantee that elements of $\mathcal{V}$ can be approximated by smooth 
infinitesimal isometries. An interesting example, discovered by Cohn-Vossen
\cite{Spivak}, is a closed smooth surface of non-negative curvature for which 
$\mathcal{C}^\infty \cap \mathcal{V}$ consists only of trivial fields 
$V: S\longrightarrow {\mathbb R}^3$ with constant gradient, 
whereas $\mathcal{C}^2 \cap \mathcal{V}$ contains non-trivial infinitesimal isometries.
Therefore $\mathcal{C}^\infty \cap \mathcal{V}$ is not dense in $\mathcal{V}$ for this surface.
We however have:
\begin{theorem}\label{th_density-intro}
Assume that $S$ is elliptic, homeomorphic to a disk, of class $\mathcal{C}^{m+2,\alpha}$ 
up to the boundary and that $\partial S$ is $\mathcal{C}^{m+1,\alpha}$, 
for some $\alpha\in (0,1)$ and an integer $m>0$. 
Then, for every $V\in\mathcal{V}$ there exists a sequence 
$V_n\in\mathcal{V}\cap\mathcal{C}^{m,\alpha}(\bar S,\mathbb{R}^3) $ 
such that: $$\lim_{n\to\infty} \|V_n - V\|_{W^{2,2}(S)} = 0.$$
\end{theorem} 
 
\medskip

Recall that the $\Gamma$-convergence result is a combination of two statements.
The first one concerns compactness and lower bound for any equi-bounded
sequence of $3$d deformations $u^h: S^h \longrightarrow{\mathbb R^3}$ 
in terms of the limit energy. The second one states an upper bound on a
recovery sequence given an admissible mapping 
for the limit problem \cite{dalmaso}. 
In view of future applications, we prove the lower bound result for an arbitrary surface. 

\begin{theorem}\label{thliminf-intro}
Let $S$ be a $2$-dimensional surface embedded in $\mathbb{R}^3$, which is
compact, connected, oriented, of class $\mathcal{C}^{1,1}$,
and whose boundary $\partial S$ is the union of finitely many 
(possibly none) Lipschitz curves. Assume (\ref{scaling-intro}) 
and let $u^h\in W^{1,2}(S^h,\mathbb{R}^3)$ be a sequence of deformations whose scaled energies 
$E^h(u^h)/e^h$ are uniformly bounded.
Then there exist a sequence $Q^h\in SO(3)$ and $c^h\in\mathbb{R}^3$ 
such that for the normalized rescaled deformations:
$$y^h(x+t\vec{n}) = Q^h u^h(x+h/h_0 t\vec{n}) - c^h$$ 
defined on the common domain $S^{h_0}$, the following holds.
\begin{enumerate}
\item[(i)] $y^h$ converge in $W^{1,2}(S^{h_0})$ to $\pi$.
\item[(ii)]  The scaled average displacements:
\begin{equation}\label{Vh-intro}
V^h(x) = \frac{h}{\sqrt{e^h}} \fint_{-h_0/2}^{h_0/2}
y^h(x+t\vec{n}) - x ~\mathrm{d}t
\end{equation} 
converge (up to a subsequence) in $W^{1,2}(S)$ to some $V\in \mathcal{V}$.
\item[(iii)] ${h}/{\sqrt{e^h}} ~\mathrm{sym }~\nabla V^h$ converge
(up to a subsequence) in $L^{2}(S)$ to ${1}/{2}(A^2)_{tan}$, where $A$ is
related to $V$ by (\ref{Adef-intro}).
Equivalently, we have: $(A^2)_{tan}\in\mathcal{B}$.
\item[(iv)] $\liminf_{h\to 0} {1}/{e^h} E^h(u^h) \geq I(V).$
\end{enumerate}
\end{theorem}
The novelty with respect to the equivalent result for $\beta =4$
\cite{lemopa1} is the constraint $(A^2)_{tan}\in\mathcal{B}$
or equivalently: $V\in {\mathcal V_2}$.
If $S$ is an elliptic surface as in (\ref{elliptic-intro}), and of sufficient
regularity, 
one can prove that the set $\mathcal B$ coincides with the whole space 
$L^2_{sym}(S, {\mathbb R}^{2\times2})$ \cite{lemopa1}, hence the constraint 
is automatically satisfied for all $V\in {\mathcal V}$.  
In the general case where $S$ is an arbitrary surface, a characterization of 
this constraint and the exact form of $\mathcal B$
may be difficult.

We conjecture that other constraints, similar to the inclusion
$V\in {\mathcal V_2}$ obtained above, should be present 
for values of $\beta$ closer to $2$ (not
derived here). One then expects such conditions to be 
necessary and sufficient for constructing the recovery
sequence on shells which are not convex.
Heuristically, the closer $\beta$ is to $2$, we expect $V$ to be an
infinitesimal isometry of higher order.

\medskip
 
We now state the upper bound in the $\Gamma$-convergence result, for elliptic surfaces:
\begin{theorem}\label{thlimsup-intro}
Under the assumptions on $S$ in Theorem \ref{thliminf-intro} and that of $S$
being elliptic, homeomorphic to a disk, 
of class $\mathcal{C}^{4,\alpha}$ up to its boundary and with $\partial S$ 
of class $\mathcal{C}^{3,\alpha}$, the following holds.
Assume \eqref{scaling-intro}. Then for every $V\in\mathcal{V}$ there exists a sequence
of deformations $u^h\in W^{1,2}(S^h,{\mathbb R}^3)$ (satisfying 
$E^h(u^h)\leq C e^h$) such that:
\begin{itemize}
\item[(i)] the rescaled deformations $y^h(x+t\vec n)=u^h(x+th/h_0 \vec n)$ 
converge in $W^{1,2}(S^{h_0})$ to $\pi$.
\item[(ii)] the scaled average displacements $V^h$ given in (\ref{Vh-intro})
converge in $W^{1,2}(S)$ to $V$.
\item[(iii)] $\lim_{h\to 0} {1}/{e^h} E^h(u^h) = I(V)$.
\end{itemize}
\end{theorem}

As an application of Theorem \ref{thliminf-intro} and Theorem
\ref{thlimsup-intro}, we will further discuss in section \ref{sec_deadloads} 
the behavior of the minimizers of the total energy functional (\ref{total-intro}).
According to the static elasticity theory, they correspond to the equilibrium
configurations of the thin shell subject to a body
force $f^h$. We identify the scaling regimes for $f^h$ under which the scaling
assumption (\ref{scaling-intro}) is satisfied for the minimizers of $J^h$ and
prove the appropriate convergence result in Theorem \ref{thmaincinque}. 

\medskip

The problem of the limiting theory in the scaling range $2<\beta<4$ is still
open for general shells. Applying methods of our present paper to surfaces
changing type leads to working with mixed-type PDEs. In a next step, we plan
to extend our study to other classes of surfaces, eg hyperbolic, surfaces of
revolution and developable surfaces. This program will hopefully
pave the way for a better understanding of the most general cases. 

\bigskip

\noindent{\bf Acknowledgments.} 
M.L. was partially supported by the NSF grant DMS-0707275. 
and by the Center for Nonlinear Analysis (CNA) under
the NSF grants 0405343 and 0635983.
M.G.M. was partially supported by MiUR through 
the project ``Variational problems with multiple scales'' 2006
and by GNAMPA through the project ``Problemi di riduzione di dimensione
per strutture elastiche sottili'' 2008. M.R.P. was partially supported 
by the University of Pittsburgh grant CRDF-9003034.

\section{A lower bound  for general surfaces: a proof of Theorem \ref{thliminf-intro}}

The claims (i), (ii) and (iv) in Theorem \ref{thliminf-intro}
were proved in \cite[Theorem 2.1]{lemopa1}, under a less
restrictive condition $\beta>2$.  It remains to now deduce (iii), relying on the assumption 
that $\lim_{h\to 0} h^2/\sqrt{e^h}=0$, which follows from (\ref{scaling-intro}).
Define: 
\begin{equation}\label{gradh}
\nabla_h y^h(x+t\vec n) = Q^h \nabla u^h(x+h/h_0 t\vec n).
\end{equation}
A crucial observation made in \cite{lemopa1} is that, 
for some sequence of matrix fields $R^h\in W^{1,2}(S, \mathbb{R}^{3\times 3})$
with values in $SO(3)$, the sequence 
$1/\sqrt{e^h} (\nabla_h y^h - Q^h R^h)$ is uniformly bounded in
$L^2(S^{h_0})$ (see \cite[Proposition 3.4 (i)]{lemopa1}).  
Hence there follows the uniform bound on 
$1/\sqrt{e^h} \fint_{-h_0/2}^{h_0/2}\nabla_h y^h - Q^h R^h~\mbox{d}t $ in
$L^2(S)$. After multiplying by $h^2/\sqrt{e^h}$, we obtain:
\begin{equation}\label{uno}
\lim_{h\to 0} \frac{h^2}{e^h}\fint_{-h_0/2}^{h_0/2}\mbox{sym } (\nabla_h y^h -
\mbox{Id})~\mbox{d}t = - \lim_{h\to 0} \frac{h^2}{e^h}\mbox{sym }(\mbox{Id} -
Q^h R^h) = \frac{1}{2} A^2 \qquad \mbox{in } L^2(S),
\end{equation}
where the last equality follows from \cite[Lemma 3.2 (iii)]{lemopa1}.

On the other hand, by \cite[Proposition 3.4 (ii)]{lemopa1}, (a subsequence of) 
$h/\sqrt{e^h}(\nabla_h y^h -\mbox{Id})$ converges  in $L^2(S^{h_0})$
(to $A\pi$), and therefore $h/\sqrt{e^h}\fint_{-h_0/2}^{h_0/2}|\nabla_h y^h -
\mbox{Id}|~\mbox{d}t$ is bounded in $L^2(S)$.
Consider the matrix fields:
$$\nabla V^h(x) = \frac{h}{\sqrt{e^h}} \fint_{-h_0/2}^{h_0/2} (\nabla_h y^h(x
+ t\vec n) -\mbox{Id}) (\mbox{Id} + h/h_0t \Pi)~\mbox{d}t $$
and notice that by the previous observation the right hand side in: 
\begin{equation*}
\frac{h}{\sqrt{e^h}}\left|\mbox{sym }\nabla V^h -
  \frac{h}{\sqrt{e^h}}\mbox{sym}_{tan}\fint_{-h_0/2}^{h_0/2}\nabla_h
  y^h-\mbox{Id}~\mbox{d}t \right|
\leq \frac{h^2}{\sqrt{e^h}}\frac{h}{\sqrt{e^h}}
\fint_{-h_0/2}^{h_0/2}|\nabla_h y^h - \mbox{Id}|~\mbox{d}t
\end{equation*}
converges to $0$ in  $L^2(S)$, again using (\ref{scaling-intro}).
Recalling (\ref{uno}) we arrive at:
\begin{equation}\label{tre}
(A^2)_{tan} = 2\lim_{h\to 0} \frac{h}{\sqrt{e^h}} \mbox{ sym }\nabla V^h 
\qquad \mbox{in } L^2(S),
\end{equation}
which ends the proof of (iii) and also establishes Theorem \ref{thliminf-intro}.
\endproof

\medskip

\begin{remark}\label{remi1}
Condition (\ref{tre}) may seem more restrictive than that of 
$(A^2)_{tan}\in\mathcal{B}$, since:
\begin{equation}\label{R1}
A_{tan} = \lim_{h\to 0} \nabla V^h\quad \mbox{ in } L^2(S).
\end{equation} 
This is however not the case, because of the following observation.
Let $A$ be any skew-symmetric matrix field on $S$, such that $A_{tan} = \nabla V$
for some $V\in W^{1,2}(S)$ and that $(A^2)_{tan}\in\mathcal{B}$.
Then there must exist a sequence $V^h\in W^{1,2}(S,\mathbb{R}^3)$
such that (\ref{R1}) and (\ref{tre}) hold.
Indeed, by ``slowing down'' the sequence $\tilde V^h$ in:
$(A^2)_{tan} = 2 \lim_{h\to 0}\mbox{ sym }\nabla \tilde V^h$ if necessary, 
we may without loss of generality assume that 
$\|\tilde V^h\|_{W^{1,2}(S)}^2\leq h/\sqrt{e^h}$.
Define $V^h= \sqrt{e^h}/h \tilde V^h + V$.
Then $h/\sqrt{e^h}\mbox{ sym }\nabla V^h = \mbox{ sym }\nabla \tilde V^h$
converges to $1/2 (A^2)_{tan}$ and $\nabla V^h = \sqrt{e^h}/h \nabla \tilde V^h 
+ A_{tan}$ converges to $A_{tan}$, as the norm of the first term, bounded by 
$(\sqrt{e^h}/h)^{1/2}$ goes to $0$ in $L^2(S)$.
\end{remark}

\begin{remark}\label{remi2}
If $S$ is a plate, then the component $V^3=V\vec n$ of $V\in\mathcal{V}$ 
and the vector field $w$ 
such that $1/2(A^2)_{tan}=\mbox{sym }\nabla w$ are, respectively, 
the out-of-plane and 
the in-plane displacements (modulo a possible in-plane infinitesimal rigid motion).
Also, the constraint in (iii) becomes: $\mbox{det}\nabla^2 V^3 = 0$.
\end{remark}

\section{The linear problem $\mbox{sym}\nabla w = B$ on 
elliptic surfaces}\label{sec_linear}

Assume that $S$ is a simply connected, compact surface of class  
$\mathcal{C}^{2,1}$ with non-empty $\mathcal{C}^2$ boundary, and that its shape operator 
$\Pi$ satisfies (\ref{elliptic-intro}).
Given a vector field $w:S\longrightarrow \mathbb{R}^3$, we will 
consider its decomposition into the tangential and normal parts:
$w=w_{tan} + (w\cdot\vec n)\vec n$. Notice that:
$$ \mbox{sym}\nabla w = \mbox{sym}\nabla w_{tan} + (w\cdot\vec n)\Pi.$$
The purpose of this section is to prove the following result:
\begin{theorem}\label{thmain_linear}
There exists a linear operator: 
$$\mathcal{T}:L^2_{sym}(S,\mathbb{R}^{2\times 2})\longrightarrow 
\big\{w\in L^2(S,\mathbb{R}^3);~ w_{tan}\in W^{1,2}(S)\big\}$$
such that $\mathrm{sym}\nabla (\mathcal{T}B) = B$, for every 
$B\in L^2_{sym}(S,\mathbb{R}^{2\times 2})$ and that:
\begin{equation}\label{main_lin_est}
\|(\mathcal{T}B)_{tan}\|_{W^{1,2}(S)} + \|(\mathcal{T}B)\cdot\vec n\|_{L^2(S)}
\leq C \|B\|_{L^2(S)}.
\end{equation}
\end{theorem}
We first notice that, by a density argument, it is clearly enough to define
 the linear operator $\mathcal{T}$ on $W^{2,2}\cap 
L^2_{sym}(S,\mathbb{R}^{2\times 2})$ 
and prove there the uniform bound (\ref{main_lin_est}).

\medskip

We will use the notation and calculations in \cite{Ni} or
\cite[Section 9.2]{HH}. Since $S$ is homeomorphic to a disk, 
it can be parameterized by a single chart  
$r\in  \mathcal{C}^{2,1} (\bar\Omega, \mathbb{R}^3)$, defined on an  
open, bounded, simply connected domain $\Omega\subset {\mathbb R}^2$ with 
$\mathcal{C}^2$ boundary ($r$ can even be a conformal parameterization by 
\cite[Theorem 3.1]{Jost-book}, 
but we do not require it). The positive definite matrix field
$[g_{ij}]\in\mathcal{C}^{1,1}(\bar\Omega,\mathbb{R}^{2\times 2})$ with  
$g_{ij} = \partial_i r \cdot \partial_j r$ is the pull-back metric on $\Omega$,  
and $\sqrt{|g|}= \sqrt {\det[g_{ij}]}\in \mathcal{C}^{1,1}(\bar\Omega)$ 
is the associated volume form.  
The shape operator $\Pi$ expressed in the flat coordinates 
is given by $[h_{ij}]\in\mathcal{C}^{0,1}(\bar\Omega, \mathbb{R}^{2\times 2})$, where 
$h_{ij} = \partial_i (\vec n\circ r)\cdot\partial_j r.$  
We denote $[h^{ij}]=[h_{ij}]^{-1}=\Pi^{-1}$ and $[g^{ij}]=[g_{ij}]^{-1}$.  
The mean curvature of $S$ is given by
$H = \frac{1}{2}\mbox{tr}([g^{ij}]\Pi) \in \mathcal{C}^{0,1}(\bar\Omega)$.   

With this notation, the problem $\mbox{sym}\nabla w = B$ for a given
$B\in W^{2,2}\cap L^2_{sym}(S,\mathbb{R}^{2\times 2})$, is equivalent to the following
system of PDEs in $\Omega$:
\begin{equation}\label{r-equation} 
\left\{\begin{array}{l}  
\partial_1 r \cdot \partial_1 w = B_{11}  \\ 
\partial_1 r \cdot \partial_2 w + \partial_2 r \cdot \partial_1 w = 2 B_{12}\\ 
\partial_2 r \cdot \partial_2 w = B_{22},  
\end{array} \right. 
\end{equation}  
where we set:
\begin{equation}\label{B}
[B_{ij}]= [\partial_i r \cdot B \partial_j r].
\end{equation}
The first step consists of studying the scalar field:
\begin{equation}\label{curl}
\omega = \frac{1}{\sqrt{|g|}}\mbox{curl } w = \frac{1}{\sqrt{|g|}} 
\left(\partial_1 w \cdot\partial_2 r  
- \partial_2 w\cdot\partial_1 r\right).
\end{equation}  
Notice that $\mbox{curl }w = \mbox{curl }w_{tan}$, as the matrix field $\Pi$
is symmetric.
Now (\ref{r-equation}) formally yields:
\begin{equation}\label{curl_eq}
\mathcal{L}\omega := -\sum_{i,j=1}^2 \partial_i(\sqrt{|g|} h^{ij}\partial_j \omega) 
- 2\sqrt{|g|} H\omega = \mathcal{D}([B_{ij}]).
\end{equation}
Recall that the coefficients of the operator $\mathcal{L}$ have the following regularity:
$h^{ij}, H\in \mathcal{C}^{0,1}(\bar\Omega)$,
$\sqrt{|g|}\in\mathcal{C}^{1,1}(\bar\Omega)$. 
The matrix field $[h^{ij}]$ is uniformly strictly positive definite and both $\sqrt{|g|}$ 
and $H$ are strictly positive in $\bar\Omega$.
The linear operator $\mathcal{D}:W^{2,2}(\Omega,\mathbb{R}^{2\times 2})
\longrightarrow L^2(\Omega)$ given explicitely in \cite{HH}, has the form:
$\mathcal{D}([B_{ij}]) = \sum \partial_k(a_{klij}\partial_l B_{ij}) 
+ \sum \partial_k(b_{kij}B_{ij}) + \sum c_{ij}B_{ij}$ with coefficients
$a,b,c\in \mathcal{C}^{0,1}(\bar\Omega)$.

\begin{theorem}\label{basic_bound}
There exists a linear operator $\mathcal{S}:
W^{2,2}(\Omega,\mathbb{R}^{2\times 2}) \longrightarrow
W^{2,2}(\Omega,\mathbb{R})$ such that:
\begin{equation}\label{ineq}
\|\mathcal{S}([B_{ij}])\|_{L^2(\Omega)}\leq C\|[B_{ij}]\|_{L^2(\Omega)},
\end{equation}
and that $\mathcal{S}([B_{ij}])$ is a solution to (\ref{curl_eq}), for each
$[B_{ij}]\in W^{2,2}(\Omega,\mathbb{R}^{2\times 2})$.
\end{theorem}
\begin{proof}
{\bf 1.} Fix a small $\epsilon>0$ and 
extend the coefficients $[h^{ij}]$, $H$, $\sqrt{|g|}$ on the domain
$\Omega_\epsilon=\{x\in\mathbb{R}^2; ~\mbox{dist}(x,\Omega)<\epsilon\}$,
preserving their regularity and sign and increasing their indicated norms 
at most by a uniform
factor $C$. Also, let $W^{2,2}(\Omega,\mathbb{R}^2)\ni [B_{ij}]\mapsto[\tilde
B_{ij}]\in W^{2,2}(\Omega_\epsilon,\mathbb{R}^2)$ be a linear extension
operator, such that:
\begin{equation}\label{ext}
\mbox{supp} [\tilde B_{ij}]\subset \bar\Omega_{\epsilon/2},
\,\,\,  \|[\tilde B_{ij}]\|_{L^2(\Omega_\epsilon)}
\leq C \|[B_{ij}]\|_{L^2(\Omega)} \quad \mbox{and}\quad
\|[\tilde B_{ij}]\|_{W^{1,2}(\Omega_\epsilon)}
\leq C \|[B_{ij}]\|_{W^{1,2}(\Omega)}.
\end{equation}
For $\lambda= 2\max_{\Omega_\epsilon}(\sqrt{|g|}H)$ consider the following
bilinear, symmetric, continuous and coercitive form $a_\lambda$ on
$W^{1,2}_0(\Omega_\epsilon)$:
\begin{equation}\label{form}
a_\lambda (\omega, v) = \int_{\Omega_\epsilon} \sum_{i,j=1}^2 \sqrt{|g|}
h^{ij}\partial_j \omega\partial_i v + \int_{\Omega_\epsilon} (\lambda - 2
\sqrt{|g|}H)\omega v.
\end{equation}
By the Lax-Milgram and Rellich-Kondrachov theorems there exists
the compact linear solution operator $T:L^2(\Omega_\epsilon)\longrightarrow
L^2(\Omega_\epsilon)$ with $Range(T)\subset W^{1,2}_0(\Omega_\epsilon)$,
defined by:
$$\forall f\in L^2(\Omega_\epsilon) \quad \forall v\in W^{1,2}_0(\Omega_\epsilon)
\qquad  a_\lambda(Tf, v) = \int_{\Omega_\epsilon}fv.$$
Further, the space $Ker(\lambda T - \mbox{Id})$ has finite dimension $d$.
Since the operator $\lambda T-\mbox{Id}$ is self-adjoint on
$L^2(\Omega_\epsilon)$, by Fredholm's alternative the closed spaces
$Ker(\lambda T - \mbox{Id})$ and $E= Range(\lambda T - \mbox{Id})$ are
orthogonally complementary in $L^2(\Omega_\epsilon)$.  Consequently, there
exists a continuous linear bijection $S:E\longrightarrow E$ which is the right
inverse of $\lambda T-\mbox{Id}$, that is:
$(\lambda T-\mbox{Id})\circ S = \mbox{Id}_{E}$.

Given $f\in E$, consider $\omega=\frac{1}{\lambda} (-S-\mbox{Id})f$. Clearly:
\begin{equation}\label{first_bd}
\|\omega\|_{L^2(\Omega_\epsilon)} \leq C \|f\|_{L^2(\Omega_\epsilon)}.
\end{equation}
Also, we have: $(\lambda T-\mbox{Id})\omega = -Tf$ which is equivalent to
$T(\lambda\omega + f) = \omega$.  
Therefore $\omega\in W^{1,2}_0(\Omega_\epsilon)$ and in view of (\ref{form}): 
$a_\lambda(\omega, v) = \int_{\Omega_\epsilon} (\lambda\omega + f)v$ 
for every $v\in W^{1,2}_0(\Omega_\epsilon)$. In other words, $\omega$ is a
weak, traceless solution to the second order elliptic PDE $\mathcal{L}\omega
= f$ in $\Omega_\epsilon$, and as such must be its classical solution 
(see \cite[Theorem 8.12]{GT}) $\omega\in W^{2,2}(\Omega_\epsilon)$. Moreover:
\begin{equation}\label{scnd_bd}
\|\omega\|_{W^{2,2}(\Omega_\epsilon)} \leq C
(\|\omega\|_{L^2(\Omega_\epsilon)} + \|f\|_{L^2(\Omega_\epsilon)}) \leq
C \|f\|_{L^2(\Omega_\epsilon)},
\end{equation}
where we have used (\ref{first_bd}) to obtain the second inequality.

\medskip

{\bf 2.} Let $Ker(\lambda T - \mbox{Id}) =
span_{L^2(\Omega_\epsilon)}\{\xi_1\ldots \xi_d\}$.  Consider the operator
$L:\mathcal{C}_c^\infty(\Omega_\epsilon \setminus
\bar\Omega_{\epsilon/2})\longrightarrow \mathbb{R}^d$ given by:
$$L(f_0) = \left\{\int_{\Omega_\epsilon}f_0 \xi_k\right\}_{k=1}^d.$$
We now prove that $L$ is surjective. For otherwise, there would be:
$\int_{\Omega_\epsilon} f_0 (\sum_{k=1}^d\alpha_k\xi_k) = 0$,
for some vector $(\alpha_1\ldots \alpha_d)\in\mathbb{R}^d\setminus \{0\}$
and all $f_0\in \mathcal{C}_c^\infty(\Omega_\epsilon \setminus
\bar\Omega_{\epsilon/2})$. This implies that $\xi= \sum_{k=1}^d\alpha_k\xi_k=0$
in $\Omega_\epsilon \setminus\bar\Omega_{\epsilon/2}$.  Moreover, as $\xi\in Ker
(\lambda T - \mbox{Id})$, one obtains that $T(\lambda\xi) = \xi$ and hence
$\mathcal{L}\xi = 0$ in $\Omega_\epsilon$.
Using H\"ormander's continuation principle we deduce that the compactly
supported $\xi$ must actually vanish on the entire $\Omega_\epsilon$,
contradicting thus the linear independence of $\{\xi_1\ldots \xi_d\}$.

In view of the above, $L$ has a right inverse $L_1:\mathbb{R}^d\longrightarrow
\mathcal{C}_c^\infty(\Omega_\epsilon \setminus \bar\Omega_{\epsilon/2})$, such
that: $L\circ L_1 = \mbox{Id}_{\mathbb{R}^d}$.

\medskip

{\bf 3.} We now set:
\begin{equation*}
\mathcal{S}([B_{ij}]) := -\frac{1}{\lambda} (S+\mbox{Id})\Bigg\{
\mathcal{D}([\tilde B_{ij}]) - L_1\left(
\left\{\int_{\Omega_\epsilon}\mathcal{D}([\tilde B_{ij}]) \xi_k\right\}_{k=1}^d
\right)\Bigg\}
\end{equation*}
(or, more precisely, the restriction of the above function to $\Omega$).
By the definition of $L_1$, the quantity 
$f= \mathcal{D}([\tilde B_{ij}]) - L_1\left(\left\{\int_{\Omega_\epsilon}
\mathcal{D}([\tilde B_{ij}]) \xi_k\right\}_{k=1}^d\right)$
is orthogonal to each $\xi_i$ and so $f\in E$.  Hence $\mathcal{S}$ is well
defined, and clearly it is also linear.

Call $\omega = \mathcal{S}([B_{ij}])$.  By step 1, $\omega\in
W^{2,2}(\Omega_\epsilon)$ and the equation $\mathcal{L}\omega = f$ holds in
$\Omega_\epsilon$.  Recalling that $Range (L_1)\subset \mathcal{C}_c^\infty
(\Omega_\epsilon \setminus \bar\Omega_{\epsilon/2})$ we deduce that
the restriction $\omega_{|\Omega}$ solves (\ref{curl_eq}) in $\Omega$.  
It now remains to prove the bound (\ref{ineq}).

Since $\omega\in E$, we may define $\omega_1 = -\frac{1}{\lambda}(S+\mbox{Id})\omega$
which is a strong traceless 
solution to $\mathcal{L}\omega_1 = \omega$, and by (\ref{scnd_bd}) it satisfies:
\begin{equation}\label{thrd_bd}
\|\omega_1\|_{W^{2,2}(\Omega_\epsilon)} \leq C \|\omega\|_{L^2(\Omega_\epsilon)}.
\end{equation}
Now integrate the equation $\mathcal{L}\omega = f$ against $\omega_1$ on
$\Omega_\epsilon$.  Integration by parts (or, in other words, the selfadjointness
of the operator $T$) yields:
\begin{equation}\label{bd_final}
\|\omega\|^2_{L^2(\Omega_\epsilon)} = \int_{\Omega_\epsilon} \omega
\mathcal{L}\omega_1 =  \int_{\Omega_\epsilon} f \omega_1 \leq
\int_{\Omega_\epsilon}\mathcal{D}([\tilde B_{ij}]) \omega_1
+ C \|f_0\|_{L^2(\Omega_\epsilon)} \|\omega_1\|_{L^2(\Omega_\epsilon)},
\end{equation}
where we denoted $f_0 = L_1\left(\left\{\int_{\Omega_\epsilon}\mathcal{D}
([\tilde B_{ij}]) \xi_k\right\}_{k=1}^d \right)$. 
We treat separately the two terms in the right hand side of
(\ref{bd_final}).
For the first one, recall that the matrix field $[\tilde B_{ij}]$  is
compactly  supported in $\Omega_\epsilon$. Hence, integrating by parts twice, we may
put both derivatives of the operator $\mathcal{D}$ on $\omega_1$ and using the
Cauchy-Schwartz inequality obtain:
$$\int_{\Omega_\epsilon}\mathcal{D}([\tilde B_{ij}]) \omega_1\leq
C \|[\tilde B_{ij}]\|_{L^2(\Omega_\epsilon)} \|\omega_1\|_{W^{2,2}(\Omega_\epsilon)}. $$
For the second term in (\ref{bd_final}), notice that:
\begin{equation}\label{ine}
\begin{split}
\|f_0\|_{L^2(\Omega_\epsilon)}&\leq \|L_1\| \left\{\sum_{k=1}^d \left(
\int_{\Omega_\epsilon}\mathcal{D}([\tilde B_{ij}])
\xi_k\right)^2\right\}^{1/2}
\leq C \left\{\sum_{k=1}^d 
\|[\tilde B_{ij}]\|^2_{L^2(\Omega_\epsilon)}\|\xi_k\|^2_{W^{2,2}(\Omega_\epsilon)}
\right\}^{1/2} \\ &\leq C \|[\tilde B_{ij}]\|_{L^2(\Omega_\epsilon)},
\end{split}
\end{equation}
where we again used integration by parts.
Hence, (\ref{bd_final}) becomes:
\begin{equation}\label{ine2}
\|\omega\|^2_{L^2(\Omega_\epsilon)}\leq C \|[\tilde B_{ij}]\|_{L^2(\Omega_\epsilon)}
\|\omega_1\|_{W^{2,2}(\Omega_\epsilon)} \leq C \|[B_{ij}]\|_{L^2(\Omega)}
\|\omega\|_{L^2(\Omega_\epsilon)},
\end{equation}
by (\ref{thrd_bd}). This clearly implies (\ref{ineq}) and ends the proof 
of Theorem \ref{basic_bound}.
\end{proof}

\begin{corollary}\label{cor3.3}
The operator $\mathcal{S}$ defined in Theorem \ref{basic_bound} satisfies
the following uniform bound:
\begin{equation*}
\|\mathcal{S}([B_{ij}])\|_{W^{1,2}(\Omega)}\leq C\|[B_{ij}]\|_{W^{1,2}(\Omega)}.
\end{equation*}
\end{corollary}
\begin{proof}
Recalling the construction and notation in the proof of Theorem \ref{basic_bound}, we have:
$\omega = \mathcal{S}([B_{ij}])\in W^{2,2}(\Omega_\epsilon) 
\cap W_0^{1,2}(\Omega_\epsilon)$ 
and $\mathcal{L}\omega = \mathcal{D}([\tilde B_{ij}]) - f_0$, where $f_0$ satisfies 
(\ref{ine}).  Integrating by parts once and using the Cauchy-Schwartz inequality, we obtain:
\begin{equation}\label{wz7}
\int_{\Omega_\epsilon}\omega\mathcal{L}\omega = 
\int_{\Omega_\epsilon}\mathcal{D}([\tilde B_{ij}])\omega - \int_{\Omega_\epsilon}f_0\omega
\leq C \|[\tilde B_{ij}]\|_{W^{1,2}(\Omega_\epsilon)}\|\omega\|_{W^{1,2}(\Omega_\epsilon)}. 
\end{equation}
On the other hand, the strict ellipticity of the leading order term in 
(\ref{curl_eq}) implies that $\int\omega\mathcal{L}\omega 
\geq C\int|\nabla\omega|^2 - {1}/{C}\int|\omega|^2$.  But by (\ref{ine2}) 
$\|\omega\|^2_{L^2(\Omega_\epsilon)}\leq C\|[\tilde B_{ij}]\|_{W^{1,2}(\Omega_\epsilon)}
\|\omega\|_{W^{1,2}(\Omega_\epsilon)}$. Hence, (\ref{wz7}) and the Poincar\'e inequality yield:
$$ \|\omega\|_{W^{1,2}(\Omega_\epsilon)}^2 \leq C \|\nabla\omega\|_{L^{2}(\Omega_\epsilon)}^2
\leq C  \|[\tilde B_{ij}]\|_{W^{1,2}(\Omega_\epsilon)}
\|\omega\|_{W^{1,2}(\Omega_\epsilon)}.$$
Consequently: $ \|\omega\|_{W^{1,2}(\Omega_\epsilon)}\leq C 
\|[\tilde B_{ij}]\|_{W^{1,2}(\Omega_\epsilon)}$, which achieves the claim in view 
of (\ref{ext}). 
\end{proof}

\medskip

Towards the proof of the main result in Theorem \ref{thmain_linear}, we will
use the following
generalized version of Korn's second inequality in $2$d. The classical Korn's
inequality \cite{KO} on surfaces \cite{Jost, LewMul} states that for every
tangent $W^{1,2}$ vector field $v$ on $S$ one has:
\begin{equation}\label{korn_clas}
\|v\|_{W^{1,2}(S)}\leq C(\|v\|_{L^2(S)} + \|\mbox{sym}\nabla v\|_{L^2(S)})
\end{equation}
Notice that (\ref{korn_clas}) may be obtained from the result below by taking
$ A\equiv \frac{\sqrt 2}{2} J$, where the matrix field $J$ is defined through 
$J\tau_1 = \tau_2$ and $J\tau_2 = -\tau_1$ for any fixed smooth orthonormal basis 
$(\tau_1, \tau_2)$ of the tangent space $T_xS$.

\begin{lemma}\label{lemma_korn}
Let $A\in \mathcal{C}^{0,1}(\bar S, \mathbb{R}^{2\times 2})$ be a tensor field
on $S$, with $\det A\neq 0$ in $\bar S$. Then there holds:
$$ \|\nabla v\|_{L^2(S)} \leq C \left(\|v\|_{L^2(S)} 
+ \|(\nabla v)_{tan} - ((\nabla v)_{tan}: A) A\|_{L^2(S)}\right)$$  
for every tangent vector field $v\in W^{1,2}(S, \mathbb{R}^3)$.
\end{lemma} 
\begin{proof}
Take $J\in \mathcal{C}^{0,1}(\bar S, \mathbb{R}^{2\times 2})$ to be any
skew-symmetric $2$-tensor field on $\bar S$, with nonvanishing determinant.
Define $\tilde v=JA^{-1}v$ and observe that:
\begin{equation}\label{gradi}
\nabla v = AJ^{-1}(\nabla\tilde v)_{tan} + \nabla(AJ^{-1})JA^{-1}v.
\end{equation}
To estimate $\nabla\tilde v$, we use (\ref{korn_clas}):
\begin{equation}\label{ktilde}
\|\nabla\tilde v\|_{L^{2}(S)}
\leq C(\|v\|_{L^2(S)} + \|\mbox{sym}\nabla \tilde v\|_{L^2(S)}).
\end{equation}
Further, $\mbox{sym}\nabla \tilde v=\mbox{sym}\left((\nabla \tilde v)_{tan} 
- ((\nabla v)_{tan}: A)J\right)$ because $J$ is skew-symmetric. Therefore:
\begin{equation*}
\begin{split}
\|\mbox{sym}\nabla\tilde v\|_{L^{2}(S)}
&\leq \|(\nabla \tilde v)_{tan} - ((\nabla v)_{tan}: A)J\|_{L^2(S)} \leq 
\|AJ^{-1}(\nabla \tilde v)_{tan} - ((\nabla v)_{tan}: A)A\|_{L^2(S)}\\
&\leq \|(\nabla v)_{tan} - ((\nabla v)_{tan}: A) A\|_{L^2(S)} + C\|v\|_{L^2(S)},
\end{split}
\end{equation*}
in view of (\ref{gradi}). Combining (\ref{gradi}), (\ref{ktilde}) and the
above estimate proves the lemma.
\end{proof}

\medskip

Let $P$ denote now the orthogonal projection of the space of $2$-tensors on
$T_xS$ (this tangent space is identified with $\mathbb{R}^{2\times 2}$), onto
$[span_{\mathbb{R}^{2\times 2}}\{\Pi(x)\}]^\perp$.  Consider the following space of
tangent vector fields:
$$\mathcal{Z} = \left\{v\in W^{1,2}(S,\mathbb{R}^3); ~ v\cdot \vec n=0
\mbox{ and } P((\nabla v)_{tan})=0\right\}. $$
Notice that for every $v\in\mathcal{Z}$ there exists a scalar field $\alpha\in
L^2(S)$ such that $(\nabla v)_{tan} = \alpha\Pi$.
\begin{lemma}\label{first_korn}
We have the following:
\begin{itemize}
\item[(i)] For every tangent vector field $v\in W^{1,2}(S,\mathbb{R}^3)$ 
there holds:
$$ \|\nabla v\|_{L^2(S)} \leq C \left(\|v\|_{L^2(S)} 
+ \|P((\nabla v)_{tan})\|_{L^2(S)}\right). $$  
\item[(ii)] The space $\mathcal{Z}$ is finite dimensional.
\item[(iii)] Let $\mathcal{P}$ be the orthogonal projection of the space
  of tangent $W^{1,2}$ vector fields, onto $\mathcal{Z}$. Then:
$$ \|v - \mathcal{P} v\|_{W^{1,2}(S)} \leq C \|P((\nabla v)_{tan})\|_{L^2(S)}, $$
for every tangent vector field $v\in W^{1,2}(S,\mathbb{R}^3)$.  
\end{itemize}
\end{lemma}
\begin{proof}
The estimate in (i) follows from Lemma \ref{lemma_korn} by taking
$A=\frac{1}{\sqrt{\Pi:\Pi}}\Pi$. 

In particular, for every $v\in \mathcal{Z}$
there holds $\|\nabla v\|_{L^2(S)} \leq C\|v\|_{L^2(S)} $. It follows that in
$\mathcal{Z}$ the $L^2$ and the $W^{1,2}$ norms are equivalent, and by 
a standard argument we obtain (ii).
For otherwise the space $(\mathcal{Z}, \|\cdot\|_{W^{1,2}(S)})$ 
would have a countable Hilbertian (orthonormal) base 
$\{e_i\}_{i=1}^\infty$ and thus necessarily
the sequence $\{e_i\}$ would converge to $0$, weakly in $W^{1,2}(S)$.
But this implies that $\lim_{h\to 0}\|e_i\|_{L^2(S)} = 0$, which
by the norms equivalence gives the same convergence in $W^{1,2}(S)$,
and a contradiction.

To deduce (iii), it is enough to prove it for $v\in \mathcal{Z}^\perp$, as the
left hand side of the desired inequality vanishes for all $v\in\mathcal{Z}$.
We now argue by contradiction. Suppose that for a sequence of tangent vector
fields $v_n\in \mathcal{Z}^\perp$ one has:
\begin{equation}\label{assu}
\|v_n\|_{W^{1,2}(S)} = 1 \quad \mbox{and} \quad 
\|P((\nabla v_n)_{tan})\|_{L^2(S)}\to 0.
\end{equation}
Then, up to a subsequence, $v_n$ converges weakly (in $W^{1,2}(S)$) to some
$v\in \mathcal{Z}^\perp$. On the other hand, by the convergence in
(\ref{assu}) there must be: $P((\nabla v)_{tan})=0$, and so
$v\in \mathcal{Z}$. Consequently $v=0$, which in turn implies that 
$\|v_n\|_{L^2(S)}$ converges to $0$.  Using (i) we now deduce that 
$\nabla v_n$  converges to $0$ strongly in $L^2(S)$. Therefore, $v_n$ converges
to $0$ strongly in  $W^{1,2}(S)$, which contradicts the first assumption in
(\ref{assu}) and ends the proof.
\end{proof}

\medskip

\noindent {\bf End of proof of Theorem \ref{thmain_linear}.}\\
{\bf 1.} From calculations in \cite[Section 9.2]{HH}, it follows that given
$B\in W^{2,2}\cap L^2_{sym}(S,\mathbb{R}^{2\times 2})$, there exists a
solution $w\in W^{2,2}(\Omega,\mathbb{R}^3)$ to (\ref{r-equation}), whose
gradient is given by the following explicit formula:
\begin{equation}\label{grad_solution}
\begin{split}
\partial_1 w &= \sum_{i,j=1}^2 g^{ij}B_{1j}\partial_i r +\frac{1}{2}
\sqrt{|g|} \mathcal{S}([B_{ij}]) \sum_{i=1}^2 g^{2i}\partial_i r  + u_1 \vec n,\\
\partial_2 w &= \sum_{i,j=1}^2 g^{ij}B_{2j}\partial_i r -\frac{1}{2}
\sqrt{|g|} \mathcal{S}([B_{ij}]) \sum_{i=1}^2 g^{1i}\partial_i r  + u_2 \vec n,
\end{split}
\end{equation}
where: 
\begin{equation}\label{us} 
u_1 = \frac{1}{2} \sqrt{|g|} \sum_{i=1}^2 h^{2i}
\left(\partial_i\mathcal{S}([B_{ij}]) - c_i\right),
\quad u_2 = -\frac{1}{2} \sqrt{|g|} \sum_{i=1}^2 h^{1j}
\left(\partial_i\mathcal{S}([B_{ij}]) - c_i\right),
\end{equation}
and the scalar fields $c_i$, $i=1,2$ are given by:
\begin{equation}\label{formci} 
c_i = \frac{1}{\sqrt{|g|}}\left(\partial_1 B_{2i} - \partial_2 B_{1i}
+ \sum_{k=1}^2 (\Gamma_{2i}^k B_{1k} - \Gamma_{1i}^k B_{2k})\right).
\end{equation}
The coefficients $\Gamma_{ij}^k$ are the Christoffel symbols, that may be calculated
from:
$$ \partial_{ij}r = \sum_{k=1}^2 \Gamma_{ij}^k\partial_k r + h_{ij} \vec n,
\qquad i,j=1..2. $$

Define the operator $\mathcal{T}_1: W^{2,2}\cap L^2_{sym}(S,
\mathbb{R}^{2\times 2}) \longrightarrow W^{2,2}(S,\mathbb{R}^3)$ 
so that $w=\mathcal{T}_1 B$ satisfies 
(\ref{grad_solution}) in $\Omega$  (with a slight abuse of notation, we do not
distinguish between $w$ as the vector field on $S$ and $w\circ r$ which is the
vector field on $\Omega$), with $[B_{ij}]$ given in (\ref{B}), and
$\int_\Omega w =0$.  
Clearly $\mathcal{T}_1$ is well defined, linear and it solves:
\begin{equation}\label{cosi1}
B = \mbox{sym} \nabla(\mathcal{T}_1B) = \mbox{sym}\nabla w_{tan} +
(w\cdot\vec n)\Pi.
\end{equation}

\medskip

{\bf 2.} We will now modify $\mathcal{T}_1$ to obtain the uniform bound
(\ref{main_lin_est}). Namely, let:
$$\mathcal{T}B = w - \mathcal{P}(w_{tan}) - \alpha\vec n,
\qquad w=\mathcal{T}_1 B,$$
where in view of $Range(\mathcal{P}) \subset \mathcal{Z}$ we set:
\begin{equation}\label{cosi}
\mbox{sym}\nabla(\mathcal{P}(w_{tan})) + \alpha\Pi =0.
\end{equation}
Consequently, $\mathcal{T}$ is linear and by (\ref{cosi}):
$\mbox{sym} \nabla(\mathcal{P}(w_{tan}) + \alpha\vec n)=0$, so 
$\mbox{sym}\nabla (\mathcal{T}B) = \mbox{sym}\nabla w = B$.

Write $(\nabla w_{tan})_{tan} =\mbox{sym}\nabla w_{tan} 
+ \mbox{skew}(\nabla w_{tan})_{tan}$. By (\ref{cosi1}) we have:
$$ \|P(\mbox{sym}\nabla w_{tan})\|_{L^2(S)}
= \|P B\|_{L^2(S)}\leq \|B\|_{L^2(S)}.$$
On the other hand, recalling that $\mathcal{S}([B_{ij}]) = \mbox{curl } w =
\mbox{curl} (w_{tan})$ given by formula (\ref{curl}), we obtain:
$$ \|\mbox{skew}(\nabla w_{tan})_{tan}\|_{L^2(S)}\leq 
C \|\mathcal{S}([B_{ij}])\|_{L^2(S)} \leq C\|B\|_{L^2(S)},$$
where we used the estimate (\ref{ineq}) of Theorem \ref{basic_bound}.
By Lemma \ref{first_korn} (iii) and the two estimates above we now deduce:
$$ \|(\mathcal{T}B)_{tan}\|_{W^{1,2}(S)} 
= \|w_{tan} - \mathcal{P}(w_{tan})\|_{W^{1,2}(S)} 
\leq C \|P(\nabla w_{tan})_{tan}\|_{L^2(S)}\leq C\|B\|_{L^2(S)}.$$
Consequently, we also have:
$$ \|(\mathcal{T}B)\cdot\vec n\|_{L^2(S)}= \|B - \mbox{sym}\nabla
(\mathcal{T}B)_{tan}\|_{L^2(S)} \leq C\|B\|_{L^2(S)},$$
which concludes the proof of (\ref{main_lin_est}).
\endproof

\medskip

\begin{corollary}\label{corolp}
For any $1<p<\infty$, there exists $C_p>0$ such that 
the following holds for the  operator $\mathcal{T}$ of Theorem \ref{thmain_linear}: 
$$ \|(\mathcal{T}B)_{tan} \|_{W^{2,p}(S)} + \|(\mathcal{T}B) \cdot
\vec n)\|_{W^{1,p}(S)} \leq C_p \|B\|_{W^{1,p}(S)}. $$ 
\end{corollary} 
 
\begin{proof} 
First we extend the whole domain and coefficients of all the equations to the domain 
$\Omega_\varepsilon$ as in the proof of Theorem \ref{basic_bound}, and
construct the linear solution operator $\mathcal{T}$ associated with the larger domain,
satisfying the bound (\ref{ineq}) there. 
We now notice that the system of first order PDEs: $\mbox{sym}\nabla v +
\alpha\Pi = B$ for the unknowns $v,\alpha$, 
is elliptic in the sense of Agmon, Douglis and Nirenberg
\cite{ADN}, as shown in \cite{GSP}.  Hence, applying the local estimate of
\cite{ADN} we conclude the result.
\end{proof} 

\section{A density result on elliptic surfaces: a proof of 
Theorem \ref{th_density-intro}}\label{sec_density}

In this section we prove the density result regarding the space $\mathcal{V}$
of infinitesimal isometries on an elliptic surface $S$.
This result is a necessary step in our analysis and will be used in section \ref{sec_recovery}.
Here, in addition to assumptions made on $S$ in section 
\ref{sec_linear}, we shall 
require that $S$ is of class $\mathcal{C}^{4,\alpha}$ up to the boundary and that
$\partial S$ is $\mathcal{C}^{3,\alpha}$, for some $\alpha\in (0,1)$.
As throughout this paper, 
we focus on the case when $S$ is homeomorphic to a disk, but the result seems to hold true 
for all elliptic surfaces modulo some technical modifications in case of surfaces with 
non-trivial topology.  

We shall prove Theorem \ref{th_density-intro} only for $m=2$ ($\mathcal{C}^{2,\alpha}$ regular 
infinitesimal isometries are dense in $\mathcal{V}$), but the same reasoning
applies to higher regularities as well.

\medskip

\noindent {\bf Proof of Theorem \ref{th_density-intro}.}

\noindent {\bf 1.} From the regularity of $S$ ($\mathcal{C}^{3,1}$ is enough for this purpose),
it follows that:
$$\mbox{sym }\nabla V_{tan} = -(V\vec n)\Pi\in W^{2,2}(S,\mathbb{R}^{2\times 2}).$$ 
Writing
$V^i = (V\tau_i)\circ r\in W^{2,2}(\Omega,\mathbb{R})$ (where $\tau_i = \partial_i r$ 
for $i=1,2$) we notice that the components of the matrix field $\mbox{sym }\nabla(V^1, V^2)$
are of the form: $\partial_{\tau_i}(V\tau_j) + \partial_{\tau_j}(V\tau_i)=
V (\partial_{\tau_i}\tau_j + \partial_{\tau_j}\tau_i)$
and hence they belong to $W^{2,2}(\Omega)$.
Recall now that second derivatives of any vector field $w:\mathbb{R}^n\longrightarrow
\mathbb{R}^n$ are linear combinations of derivatives of its symmetric gradient:
$[\nabla^2 w^k]_{ij} = \partial_i [\mbox{sym }\nabla w)]_{kj} + 
\partial_j [\mbox{sym }\nabla w]_{ki} - \partial_k [\mbox{sym }\nabla w]_{ij}$.
Therefore $V^1, V^2\in W^{3,2}(\Omega)$, and so $V_{tan}\in W^{3,2}(S,\mathbb{R}^3)$.

Following calculations in section \ref{sec_linear}, we will consider the scalar 
field as in (\ref{curl}):
\begin{equation}\label{wz2}
\omega=\frac{1}{\sqrt{|g|}}\mbox{curl }V = \frac{1}{\sqrt{|g|}}\mbox{curl }V_{tan}
\in W^{2,2}(S,\mathbb{R}), 
\end{equation}
satisfying $\mathcal{L}\omega=0$. The operator $\mathcal{L}$ is defined in 
(\ref{curl_eq}) and its coefficients have regularity:
$h^{ij}, H\in \mathcal{C}^{2,\alpha}(\bar\Omega)$, 
$\sqrt{|g|}\in\mathcal{C}^{3,\alpha}(\bar\Omega)$.  Our first goal is to approximate
$\omega$ in $W^{2,2}(\Omega)$ by a sequence 
$\omega_n\in\mathcal{C}^{2,\alpha}(\bar\Omega)$ such that $\mathcal{L}\omega_n = 0$.

\medskip

{\bf 2.} Consider the space: 
$$F=\left\{v\in W^{2,2}(\Omega); ~ \mathcal{L}v = 0 \mbox{ in } \Omega,
~ v=0 \mbox{ on } \partial\Omega \right\}.$$
By \cite[Theorem 9.19]{GT}, the assumed regularity of $S$ and $\partial S$
guarantee that $F\subset \mathcal{C}^{3,\alpha}(\bar\Omega)$.

Also, $F$ has finite dimension. One easy way of seeing it is by integrating 
the equation $\mathcal{L}v=0$ against $v$ on $\Omega$. In view of the ellipticity of
$[h^{ij}]$ we obtain: $\|\nabla v\|_{L^2}\leq C \|v\|_{L^2}$, for every $v\in F$.
Therefore in $F$ the $L^2$ and the $W^{1,2}$ norms are equivalent, proving 
the claim.

Define now the finite dimensional space of traces:
$$F_{tr} = \left\{\left(\sum_{i,j=1}^2\sqrt{|g|}h^{ij}n^i\partial_j v
\right)_{\mid\partial\Omega}; ~ v\in F \right\}\subset L^2(\partial\Omega),$$
where by $n=(n^1, n^2)$ we denote the outer unit normal to $\partial\Omega$.
As a consequence of the regularity of vector fields in $F$ and the regularity of $S$,
there holds $F_{tr}\subset \mathcal{C}^{2,\alpha}(\partial\Omega)$.

The significance of the space $F_{tr}$ is the following. Given $w\in W^{2,2}(\Omega)$,
the problem: 
$$\mathcal{L}u = \mathcal{L}w \mbox{ in }\Omega, \qquad u=0 \mbox{ on }\partial\Omega$$
has a solution if and only if $w_{\mid\partial\Omega}\in F_{tr}^\perp$. 
Indeed, by Fredholm's alternative there must be $\mathcal{L}w\in F^\perp$ 
(where the orthogonal complement is taken in $L^2(\Omega)$), that is:
\begin{equation}\label{ms}
\forall v\in F \qquad \int_\Omega \mathcal{L}w\cdot v = 0.
\end{equation}
Integrating by parts we obtain:
$$\int_\Omega \mathcal{L}w\cdot v = \int_\Omega w\cdot \mathcal{L}w
+ \int_{\partial\Omega} \left(\sum_{i,j=1}^2\sqrt{|g|}h^{ij}n^j\partial_i v\right) w
= \int_{\partial\Omega} \left(\sum_{i,j=1}^2\sqrt{|g|}h^{ij}n^i\partial_j v\right) w,$$
and thus (\ref{ms}) is equivalent to: $w_{\mid\partial\Omega}\in F_{tr}^\perp$.

\medskip

{\bf 3.} Recall that $\omega\in W^{2,2}(\Omega)$, given in (\ref{wz2}), satisfies 
$\mathcal{L}\omega = 0$.  Hence the continuous function 
$\phi=\omega_{\mid\partial\Omega}$ belongs to $F_{tr}^\perp$.  Towards approximating
$\omega$ by $\mathcal{C}^{2,\alpha}(\bar\Omega)$ functions, we first approximate 
its trace $\phi$.  Namely, we claim that there exists a sequence 
$\phi_n \in F_{tr}^\perp\cap \mathcal{C}^{2,\alpha}(\partial\Omega)$ such that:
\begin{equation}\label{wz2.5}
\lim_{n\to\infty} \,\inf\left\{\|\Phi\|_{W^{2,2}(\Omega)}; ~ \Phi\in W^{2,2}(\Omega),
~ \Phi_{\mid\partial\Omega}=\phi_n-\phi\right\}=0,
\end{equation}
(it is clear that ``$\inf$'' above may be replaced by ``$\min$'').
Let $\tilde\phi_n=\tilde\omega_{n}{}_{\mid\partial\Omega}$, so that 
$\tilde\omega_n\in\mathcal{C}^{2,\alpha}(\bar\Omega)$ and:
\begin{equation}\label{wz3}
\lim_{n\to\infty} \|\tilde\omega_n - \omega\|_{W^{2,2}(\Omega)} = 0.
\end{equation}
Define $\phi_n = (\mbox{Id} - P_{tr})(\tilde\phi_n)$ where 
$P_{tr}:L^2(\partial\Omega)\longrightarrow F_{tr}$ is the orthogonal projection onto
$F_{tr}$.  Since $\lim_{n\to\infty}\tilde\phi_n = \phi$ in $L^2(\partial\Omega)$, then
$\tilde\phi_n - \phi_n = P_{tr}(\tilde\phi_n)$ converges in $L^2(\partial\Omega)$
to $P_{tr}(\phi)=0$.

In a finitely dimensional space $F_{tr}$ all norms are equivalent, hence in particular:
$$\lim_{n\to\infty} \,\inf\left\{\|\Phi\|_{W^{2,2}(\Omega)}; ~ \Phi\in W^{2,2}(\Omega),
~ \Phi_{\mid\partial\Omega}=\tilde\phi_n-\phi_n\right\}=0.$$
Together with (\ref{wz3}) the above implies (\ref{wz2.5}) and proves the claim.

\medskip

{\bf 4.} Consider now the sequence of harmonic functions $\Phi_n\in W^{2,2}(\Omega)$
such that $\Delta\Phi_n=0$ in $\Omega$, $\Phi_n = \phi_n-\phi$ on $\partial\Omega$.
By (\ref{wz2.5}) and usual elliptic estimates it follows that:
\begin{equation}\label{wz4}
\lim_{n\to\infty} \|\Phi_n\|_{W^{2,2}(\Omega)} = 0.
\end{equation}
Since $\phi_n-\phi\in F_{tr}^\perp$, the problem:
$\mathcal{L}u_n = -\mathcal{L}\Phi_n$ in $\Omega$, $u_n=0$ on $\partial\Omega$,
has a solution whose regularity, in view of \cite[Theorem 8.12]{GT} must be
$u_n\in W^{2,2}(\Omega)$. Moreover, using (\ref{wz4}):
\begin{equation}\label{wz5}
\|u_n\|_{W^{2,2}} \leq C\|\mathcal{L}\Phi_n\|_{L^2}
\leq C\|\Phi_n\|_{W^{2,2}}\to 0 \quad \mbox{ as } n\to \infty.
\end{equation}
Set $\omega_n = u_n + \omega +\Phi_n$. By (\ref{wz4}) and (\ref{wz5}):
\begin{equation}\label{wz6}
\lim_{n\to\infty}\|\omega_n-\omega\|_{W^{2,2}(\Omega)} =0.
\end{equation}
On the other hand $\mathcal{L}\omega_n=\mathcal{L}\omega = 0$ and 
${\omega_n}_{\mid\partial\Omega} = \phi_n \in\mathcal{C}^{2,\alpha}(\partial\Omega)$
so by \cite[Theorem 9.19]{GT} there must be 
$\omega_n\in\mathcal{C}^{2,\alpha}(\bar\Omega)$. 

\medskip

{\bf 5.} Based on the approximation $\{\omega_n\}$ accomplished in the previous step,
we now construct the desired approximation 
$V_n\in\mathcal{V}\cap \mathcal{C}^{2,\alpha}(\Omega)$ of $V$.
According to formulas (\ref{grad_solution}) and (\ref{us})
(see formulas (9.2.17) and (9.2.12) in \cite{HH}), where we put $[B_{ij}]=0$ to have
$\mbox{sym }\nabla V_n=0$, and $\mathcal{S}([B_{ij}])=\omega_n$ in view of 
$\mathcal{L}\omega_n=0$, the gradients $\nabla V_n$
are given by:
\begin{equation*}
\partial_1 V_n  = \frac{1}{2}\omega_n\sqrt{|g|} \sum_{i=1}^2 g^{2i}\partial_i r
+ u_1\vec n,\qquad
\partial_2 V_n  = -\frac{1}{2}\omega_n\sqrt{|g|} \sum_{i=1}^2 g^{1i}\partial_i r
+ u_2\vec n,
\end{equation*}
with: 
$$u_1= \frac{1}{2}\sqrt{|g|} \sum_{i=1}^2 h^{2i}\partial_i \omega_n, \qquad
u_2= -\frac{1}{2}\sqrt{|g|} \sum_{i=1}^2 h^{1i}\partial_i \omega_n. $$
We see that $\nabla V_n\in\mathcal{C}^{1,\alpha}(\bar\Omega)$ and 
$\lim_{n\to\infty}\|\nabla V_n - \nabla V\|_{W^{1,2}} = 0$ by (\ref{wz6}).
Normalizing $V_n$ so that $\int_{\Omega} V_n = \int_{\Omega} V$, the theorem
follows in view of Poincar\'e's inequality.
\endproof

\medskip

\begin{remark}\label{remi3}
We also expect that the assumption on the required $\mathcal{C}^{4,\alpha}$ 
regularity of $S$ could be relaxed. A less regular approximating sequence 
could be perhaps used in carrying out the analysis in sections  \ref{sec_exact} 
and \ref{sec_recovery}.  
\end{remark}

\section{Matching infinitesimal to exact isometries: 
a proof of Theorem \ref{th_exact-intro}}\label{sec_exact}

The applicability of Theorem \ref{th_exact-intro} in our analysis can be
understood by considering the limit theory of nonlinear 
elasticity for thin shells, and the scaling factor $\beta=2$. As shown in \cite{FJMM_cr}. 
one obtains then the pure bending theory of Kirchhoff, postulating that the set of 
admissible deformations is that of isometric immersions of $S$. 
Heuristically, the closer $\beta$ is to $2$, the closer  
a recovery sequence for the proposed theory should be to an exact isometry of $S$.
An evident strategy is hence to construct an exact isometry from the data of the 
limiting problem, for the whole range of $2<\beta<4$. Indeed, for elliptic shells,
such construction reduces the derived theory to 
the linear one, as explained in our paper. On the other hand, the plate 
theory for the same scaling regime \cite{FJMhier}, postulates that 
(only) the second order infinitesimal isometries of a flat domain can be modified 
(up to a second order change) into an exact isometry. This implies a
qualitatively different, nonlinear response of plates in this regime.  

\medskip

The proof of Theorem \ref{th_exact-intro} 
is based on a fixed point argument and is inspired by the proof of 
openness of the set of positive 
curvature  metrics $g$ on the topological sphere $\Sigma$ which are pull-backs
of the Euclidean metric by 
immersions of $\Sigma$ in ${\mathbb R}^3$ \cite{Ni}. 
To draw a parallel here, given an immersion 
$r$ of the topological disk with pull-back metric $g$, we seek a new immersion with the 
same pull-back metric whose first order difference from $r$ is a given 
infinitesimal isometry. 

\medskip

First, we prove a H\"older estimate for the operator $\mathcal{T}$ defined 
in section \ref{sec_linear}.

\begin{lemma}\label{Thold}
One has the following uniform estimate:
\begin{equation}\label{Tproduct}
\forall \phi,\psi\in\mathcal{C}^{2,\alpha}(\bar S,\mathbb{R}^3)
\qquad \left\|\mathcal{T}\big(\mathrm{sym}~((\nabla\phi)^T\nabla\psi)\big)
\right\|_{C^{2,\alpha}(\bar S)} \leq C \|\phi\|_{C^{2,\alpha}(\bar S)}  
\|\psi\|_{C^{2,\alpha}(\bar S)}.
\end{equation}
\end{lemma}

\begin{proof}
{\bf 1.} Call $w= \mathcal{T}(B)$ and $\omega=\mathcal{S}([B_{ij}])$, 
where $[B_{ij}]$ is given as in (\ref{B})
for $B=\mathrm{sym}~((\nabla\phi)^T\nabla\psi)$. Similarly as in Corollary
\ref{corolp}, by the interior estimates in \cite{ADN} we obtain:
\begin{equation}\label{wz8}
\|w\|_{\mathcal{C}^{1,\alpha}(\bar S)} \leq 
C\left(\|B\|_{\mathcal{C}^{1,\alpha}(\bar S)} + \|w\|_{L^\infty(S)}\right) \leq
C \|B\|_{\mathcal{C}^{1,\alpha}(\bar S)},
\end{equation}
where we used Corollary \ref{corolp} for deducing:
$\|w\|_{L^\infty}\leq \|w\|_{W^{1,p}}\leq C \|B\|_{W^{1,p}} \leq 
\|B\|_{\mathcal{C}^{1,\alpha}}$, with $p>2$.
In the same manner, the Schauder estimates for elliptic systems in divergence form \cite{GT}
Theorem 8.32 and Corollary \ref{cor3.3} also yield:
\begin{equation}\label{wz9}
\begin{split}
\|\omega\|_{\mathcal{C}^{1,\alpha}(\bar\Omega)} \leq C\left( 
\|[B_{ij}]\|_{\mathcal{C}^{1,\alpha}(\bar\Omega)} 
+ \|\omega\|_{L^\infty(\Omega)}\right)
\leq 
C \|[B_{ij}]\|_{\mathcal{C}^{1,\alpha}(\bar\Omega)}. 
\end{split}
\end{equation}
Hence, by (\ref{grad_solution}) and (\ref{wz8}):
\begin{equation}\label{wz10}
\|w\|_{\mathcal{C}^{2,\alpha}(\bar S)} \leq 
C\left(\|[B_{ij}]\|_{\mathcal{C}^{1,\alpha}(\bar \Omega)} +
\|u_1\|_{\mathcal{C}^{1,\alpha}(\bar \Omega)} +
\|u_2\|_{\mathcal{C}^{1,\alpha}(\bar \Omega)}\right).
\end{equation}
Clearly, $\|[B_{ij}]\|_{\mathcal{C}^{1,\alpha}(\bar \Omega)}\leq
C\|B\|_{\mathcal{C}^{1,\alpha}(\bar S)}$ is bounded by the right hand side of
(\ref{Tproduct}). It remains therefore to derive suitable bounds 
on the correction coefficients $u_1$ and $u_2$.

\medskip

{\bf 2.} Directly from (\ref{us}) and (\ref{formci}), in view of (\ref{wz9}) we get:
\begin{equation}\label{wz11}
\|(u_1,u_2)\|_{\mathcal{C}^{0,\alpha}(\bar \Omega)} \leq
C\left(\|\omega\|_{\mathcal{C}^{1,\alpha}(\bar \Omega)} +
\|[B_{ij}]\|_{\mathcal{C}^{1,\alpha}(\bar \Omega)}\right)
\leq C\|[B_{ij}]\|_{\mathcal{C}^{1,\alpha}(\bar \Omega)}.
\end{equation}
Also, one can check that $u_1$ and $u_2$ solve the following first order elliptic 
system (see calculations in \cite{HH} leading to (9.2.18) and (9.2.28)):
\begin{equation*}
\begin{split}
& \partial_1 u_2 - \partial_2 u_1 = \sqrt{|g|} H \omega +\frac{1}{2}\sum_{i,j=1}^2
g^{ij}(h_{j2}B_{1i} - h_{j1}B_{2i})\\
&\sum_{i,j=1}^2 \partial_i\left(\frac{\mbox{det }\Pi}{\sqrt{|g|}}h^{ij}u_j\right)
= \frac{1}{2}(\partial_2 c_1 - \partial_1 c_2).
\end{split}
\end{equation*}
As in the proof of Corollary \ref{corolp}, the interior estimates in \cite{ADN}
yield:
\begin{equation}\label{wz12}
\begin{split}
\|(u_1, u_2)\|_{\mathcal{C}^{1,\alpha}(\bar \Omega)}&\leq
C\left(\|(u_1,u_2)\|_{L^\infty(\Omega)} +
\|\omega\|_{\mathcal{C}^{0,\alpha}(\bar \Omega)} +
\|[B_{ij}]\|_{\mathcal{C}^{0,\alpha}(\bar \Omega)}  
+ \|\partial_2 c_1 - \partial_1 c_2\|_{\mathcal{C}^{0,\alpha}(\bar \Omega)}\right)\\
&\leq C\left(\|[B_{ij}]\|_{\mathcal{C}^{0,\alpha}(\bar \Omega)}  
+ \|\partial_2 c_1 - \partial_1 c_2\|_{\mathcal{C}^{0,\alpha}(\bar \Omega)}\right),
\end{split}
\end{equation}
the last bound being a consequence of (\ref{wz9}) and (\ref{wz11}).

By a direct calculation from (\ref{formci}) we obtain:
\begin{equation}\label{wz13}
\|\partial_2 c_1 - \partial_1 c_2\|_{\mathcal{C}^{0,\alpha}(\bar \Omega)} \leq
C\left(\|[B_{ij}]\|_{\mathcal{C}^{1,\alpha}(\bar \Omega)}  
+ \|\partial_{11} B_{22}+\partial_{22} B_{11} 
- 2\partial_{12}B_{22}\|_{\mathcal{C}^{0,\alpha}(\bar \Omega)}\right).
\end{equation}
Recall that: $B_{ij} =\partial_i r\cdot B\partial_j r = \frac{1}{2} (\partial_i r
\cdot (\nabla\phi)^T \nabla\psi\partial_j r +\partial_i r
\cdot (\nabla\psi)^T \nabla\phi\partial_j r) 
= \frac{1}{2} (\partial_i\phi\cdot\partial_j\psi + \partial_j\phi\cdot\partial_i\psi)$,
where, with a slight abuse of notation, we identify $\phi$ and $\psi$ with
vector fields $\phi\circ r$ and $\phi\circ r$ defined on $\bar \Omega$.
The quantity in the last term in (\ref{wz13}) equals to:
\begin{equation*}
\begin{split}
\partial_{11} B_{22}+\partial_{22} B_{11} - 2\partial_{12}B_{22}
& = \partial_{11} (\partial_2\phi\cdot\partial_2\psi) +
\partial_{22} (\partial_1\phi\cdot\partial_1\psi)
- \partial_{12} (\partial_1\phi\cdot\partial_2\psi 
+ \partial_2\phi\cdot\partial_1\psi)\\
&= - (\partial_{11}\phi\cdot\partial_{22}\psi + \partial_{22}\phi\cdot\partial_{11}\psi
- 2 \partial_{12}\phi\cdot\partial_{12}\psi).
\end{split}
\end{equation*}
Thus, by (\ref{wz13}) and (\ref{wz12}) it follows that:
\begin{equation*}
\|(u_1, u_2)\|_{\mathcal{C}^{1,\alpha}(\bar \Omega)}\leq
C\left(\|[B_{ij}]\|_{\mathcal{C}^{1,\alpha}(\bar \Omega)}  
+ \|\phi\|_{C^{2,\alpha}(\bar \Omega)} \|\psi\|_{C^{2,\alpha}(\bar \Omega)}\right)
\leq C \|\phi\|_{C^{2,\alpha}(\bar S)} \|\psi\|_{C^{2,\alpha}(\bar S)},
\end{equation*}
which completes the proof in view of(\ref{wz10}).
\end{proof}

\medskip

\noindent {\bf Proof of Theorem \ref{th_exact-intro}.}\\
For vector field $u_h = \mbox{id} + hV + h^2 w_h\in\mathcal{C}^{2,\alpha}
(\bar{S},\mathbb{R}^3)$ to be an isometry, the following condition must be satisfied:
$|\partial_\tau u_h|^2 = 1$ for every $\tau\in T_x S$.
Since:
$$ |\partial_\tau u_h|^2 - 1 = 
|\tau + h\partial_\tau V + h^2\partial_\tau w_h|^2 - 1 = h^2|\partial_\tau V|^2
+ 2h^3 \partial_\tau V\cdot \partial_\tau w_h + h^4 |\partial_\tau w_h|^2
+ 2h^2 \tau\cdot \partial_\tau w_h,$$
$u_h$ is therefore an exact isometry if and only if:
\begin{equation}\label{m1}
\mbox{sym }\nabla w_h = -\frac{1}{2} \left(\nabla V + h \nabla w_h\right)^T 
\left(\nabla V + h \nabla w_h\right).
\end{equation}
Consider the mapping $\mathcal{G}_h:\mathcal{C}^{2,\alpha}(\bar S,\mathbb{R}^3) 
\longrightarrow \mathcal{C}^{2,\alpha}(\bar S,\mathbb{R}^3)$ defined by:
$$\mathcal{G}_h(w) = -\frac{1}{2}\mathcal{T}\big(
\left(\nabla V + h \nabla w\right)^T 
\left(\nabla V + h \nabla w\right)\big). $$
By Lemma \ref{Thold} one has: $\|\mathcal{G}_h(w)\|_{\mathcal{C}^{2,\alpha}}
\leq C\|V + hw\|_{\mathcal{C}^{2,\alpha}}^2$.  Hence, putting 
$R=\|V\|_{\mathcal{C}^{2,\alpha}}^2 +1$, and denoting $\bar B_R$ the closed ball 
of radius $R$ in $\mathcal{C}^{2,\alpha}(\bar S)$,
it follows that, for all $h$ small enough, $\mathcal{G}_h(\bar B_R) \subset \bar B_R$.
We will show that $\mathcal{G}_h$ is a contraction on $\bar B_R$.
By the linearity of $\mathcal{T}$ there must be:
\begin{equation*}
\begin{split}
\mathcal{G}_h(w) - \mathcal{G}_h(v) & = -\frac{1}{2}\mathcal{T}\big(
\left(\nabla V + h\nabla w\right)^T\left(\nabla V + h\nabla w\right)
- \left(\nabla V + h\nabla v\right)^T\left(\nabla V + h\nabla v\right)\big)\\
& = -\frac{1}{2}\mathcal{T}\big( 2h \mbox{sym } (\nabla V)^T(\nabla w - \nabla v)
+ h^2 (\nabla w)^T (\nabla w - \nabla v) 
+ h^2 (\nabla w - \nabla v)^T \nabla v \big).
\end{split}
\end{equation*}
The matrix field in the argument of $\mathcal{T}$ above is clearly symmetric and 
it has the form allowing for use of Lemma \ref{Thold}. Thus, for every
$w,v\in\bar B_R$ there holds:
\begin{equation*}
\begin{split}
\|\mathcal{G}_h(w) - \mathcal{G}_h(v)\|_{\mathcal{C}^{2,\alpha}(\bar S)} &\leq
C\Big(h\|V\|_{\mathcal{C}^{2,\alpha}} + h^2 \|w\|_{\mathcal{C}^{2,\alpha}}
+ h^2 \|v\|_{\mathcal{C}^{2,\alpha}}\Big) 
\|w - v\|_{\mathcal{C}^{2,\alpha}(\bar S)}\\
&\leq Ch\big(\|V\|_{\mathcal{C}^{2,\alpha}} + 2hR\big) 
\|w - v\|_{\mathcal{C}^{2,\alpha}(\bar S)}.
\end{split}
\end{equation*}
By the Banach fixed point theorem we now conclude that for all small $h$ the problem
(\ref{m1}) has a solution $w_h$, such that $\|w_h\|_{\mathcal{C}^{2,\alpha}(\bar S)}
\leq R$. 
\endproof

\section{Construction of the recovery sequence: a proof of Theorem 
\ref{thlimsup-intro}}\label{sec_recovery}

Here, we establish the limsup part of our $\Gamma$-convergence result, 
through constructing a recovery sequence for elliptic surfaces, based on 
Theorem \ref{th_density-intro} and Theorem \ref{th_exact-intro}.

\medskip

\noindent {\bf Proof of Theorem \ref{thlimsup-intro}}

\noindent By the density result proved in Theorem \ref{th_density-intro} and the continuity of 
the functional $I$ with respect to the strong topology of $W^{2,2}(S)$, we can assume 
$V\in {\mathcal V}\cap{\mathcal C}^{2,\alpha}(\bar S, {\mathbb R}^3)$. 
In the general case the result will then follow from a diagonal argument.

In the sequel, by the Landau symbols $\mathcal{O}(s)$ and $o(s)$ we shall
denote, respectively, any quantity whose quotient with $s$ is uniformly 
bounded or converges to $0$, as $s\to 0$.

{\bf 1.} 
Let $\e=\sqrt{e^h}/h$. We recall that $\e\to0$ as $h\to0$, by assumption \eqref{scaling-intro}.
Therefore, by Theorem \ref{th_exact-intro} there exists a sequence
$w_\e:\bar S\longrightarrow \mathbb{R}^3$, equibounded in 
$\mathcal{C}^{2,\alpha}(\bar S)$, such that for all small $h>0$ the map:
\begin{equation}\label{ex-iso}
u_\e = \mathrm{id} +\e V + \e^2w_\e 
\end{equation}
is an exact isometry. 

For every $x\in S$, let $\vec n_\e(x)$ denote the unit normal vector
to $u_\e(S)$ at the point $u_\e(x)$. By the regularity of $u_\e$ we have that
$\vec n_\e\in {\mathcal C}^{1,\alpha}(\bar S,{\mathbb R}^3)$, while by \eqref{ex-iso} 
we obtain the expansion:
\begin{equation}\label{ne-exp}
\vec n_\e = \vec n +\e A\vec n + \mathcal{O}(\e^2). 
\end{equation}
Indeed, one can take $\vec n_\e = \partial_{\tau_1}u_\e \times \partial_{\tau_2}u_\e$,
where $\tau_1,\tau_2\in T_xS$ are such that $\vec n=\tau_1\times\tau_2$. Using
now the Jacobi identity for vector product and the fact that $A\in so(3)$, we
arrive at (\ref{ne-exp}).

Here we introduce the recovery sequence $u^h$ as required by the statement of the theorem. 
Note that the following suggestion for $u^h$ is in accordance with the one
used in \cite{FJMM_cr} in the framework of the purely nonlinear bending theory
for shells, corresponding to the scaling regime 
$\beta=2$. Consider the sequence of deformations $u^h\in W^{1,2}(S^h,{\mathbb R}^3)$ defined by:
\begin{equation}\label{rec_seq}
u^h(x+t\vec n) = u_\varepsilon(x) + t\vec n_\varepsilon(x)
+ \frac{t^2}{2}\varepsilon d^h(x).
\end{equation}
The vector field $d^h\in W^{1,\infty}(S,\mathbb{R}^3)$ is defined so that:
\begin{equation}\label{nd01h}
\lim_{h\to 0} h^{1/2} \|d^h\|_{W^{1,\infty}(S)} = 0,
\end{equation}
and:
\begin{equation}\label{warp}
\lim_{h\to 0} d^h = 2c\left(x,{\rm sym}(\nabla(A\vec n) - A\Pi)_{tan}\right) \quad
\mbox{ in } L^\infty(S),
\end{equation}
where $c(x,F_{tan})$ denotes the unique vector satisfying 
${\mathcal Q}_2(x, F_{tan})={\mathcal Q}_3(F_{tan} +c\otimes \vec n(x)+\vec
n(x)\otimes c)$ (see \cite[Section~6]{lemopa1}). We observe that, as 
$V\in{\mathcal C}^{2,\alpha}(\bar S, {\mathbb R}^3)$ and $c$ depends linearly 
on its second argument, the vector field:
\begin{equation}\label{def-zeta}
\zeta(x)=c(x, {\rm sym}(\nabla(A\vec n) - A\Pi)_{tan})
\end{equation}
belongs to $L^\infty(S,\mathbb{R}^3)$.

Properties (i) and (ii) now easily follow from the uniform bound on $w_\e$ 
and the normalization \eqref{nd01h}.

{\bf 2.} 
To prove (iii) it is convenient to perform a change of variables in 
the energy $E^h(u^h)$, so to express it in terms of the scaled deformation $y^h$. 
By a straightforward calculation, we obtain:
\begin{equation}\label{int-nv}
\frac{1}{e^h} E^h(u^h) = \frac{1}{e^h}\int_S \fint_{-h_0/2}^{h_0/2} 
W(\nabla_h y^h(x+t\vec n))\det[\mbox{Id}+th/h_0\Pi(x)]~\mbox{d}t\mbox{d}x,
\end{equation}
where $\nabla_h y^h(x+t\vec n)=\nabla u^h(x+th/h_0\vec n)$, as in (\ref{gradh}).
We also have:
\begin{equation}\label{form2}
\begin{split}
\nabla_h y^h(x + t\vec n) \vec n(x) & = \frac{h_0}{h} \partial_{\vec n}y^h(x+t\vec n)
= \vec n_\e(x) + th/h_0\e d^h(x),\\
\nabla_h y^h(x + t\vec n)\tau & = \nabla y^h(x+t\vec n)\cdot (\mbox{Id} +
t\Pi(x)) (\mbox{Id} + th/h_0\Pi(x))^{-1}\tau \\
& = \Big(\nabla u_\e(x) +th/h_0\nabla\vec n_\e(x) + \frac{t^2}{2h_0^2}h^2\e
\nabla d^h(x)\Big)(\mbox{Id} + th/h_0\Pi(x))^{-1}\tau, 
\end{split}
\end{equation}
for all $x\in S$ and $\tau\in T_xS$.

From \eqref{ex-iso}, \eqref{ne-exp} and \eqref{nd01h} it follows that 
$\|\nabla_h y^h-\mbox{Id}\|_{L^\infty(S^{h_0})}$ converges to $0$ as $h\to 0$. 
It now follows by polar decomposition theorem (for $h$ sufficiently small)
that $\nabla_h y^h$ is a product of a proper rotation and the well defined
square root of $(\nabla_h y^h)^T\nabla_h y^h$. By frame indifference of $W$ we
deduce that:
$$ W(\nabla_h y^h) = W\left(\sqrt{(\nabla_h y^h)^T\nabla_h y^h}\right)
= W\left(\mbox{Id} + \frac{1}{2} K^h + \mathcal{O}(|K^h|^2)\right), $$ 
where the last equality follows by Taylor expansion, with $K^h$ given by
$$ K^h =  (\nabla_h y^h)^T\nabla_h y^h - \mbox{Id}. $$
As $\|K^h\|_{L^\infty(S^{h_0})}$ is infinitesimal as $h\to 0$, we can expand
$W$ around $\mbox{Id}$ and obtain:
\begin{equation}\label{Wdopp}
\frac{1}{e^h} W(\nabla_h y^h) = \frac{1}{2} \mathcal{Q}_3\left(\frac{1}{2\sqrt{e^h}} K^h + 
\frac{1}{\sqrt{e^h}}\mathcal{O}(|K^h|^2)\right) + 
\frac{1}{\sqrt{e^h}}\mathcal{O}(|K^h|^2).
\end{equation}
Using \eqref{form2} we now calculate $K^h$. We first consider the tangential
minor of $K^h$:
\begin{equation*}
\begin{split}
K^h_{tan}(x + t\vec n) & =  (\mbox{Id} + th/h_0\Pi)^{-1}\Big[\mbox{Id}
+2 th/h_0\, {\rm sym}((\nabla u_\e)^T\nabla\vec n_\e)\\
& \qquad + t^2 h^2/h_0^2(\nabla \vec n_\e)^T\nabla\vec n_\e
+ o(\sqrt{e^h})\Big](\mbox{Id} + th/h_0\Pi)^{-1} - \mbox{Id}\\
& = (\mbox{Id} + th/h_0\Pi)^{-1}\Big[
2 th/h_0\, {\rm sym}((\nabla u_\e)^T\nabla\vec n_\e) -2 th/h_0\Pi\\
& \qquad + t^2 h^2/h_0^2(\nabla \vec n_\e)^T\nabla\vec n_\e - t^2 h^2/h_0^2\Pi^2
\Big](\mbox{Id} + th/h_0\Pi)^{-1} + o(\sqrt{e^h}),
\end{split}
\end{equation*}
where we have used the fact that $u_\e$ is an isometry to see that 
$(\nabla u_\e)^T\nabla u_\e = \mbox{Id}$, and the identity:
$$F_1^{-1}FF_1^{-1}-\mbox{Id} = F_1^{-1}(F-F_1^2)F_1^{-1}.$$ 
By \eqref{ex-iso} and \eqref{ne-exp} we also deduce:
\begin{equation*}
\begin{split}
{\rm sym} ((\nabla u_\e)^T\nabla\vec n_\e) & 
= \Pi +\e\, {\rm sym}(\nabla(A\vec n)-A\Pi) +{\mathcal O}(\e^2), \\
(\nabla \vec n_\e)^T\nabla\vec n_\e & = \Pi^2+{\mathcal O}(\e).
\end{split}
\end{equation*}
Combining these two identities with the expression of $K^h_{tan}$ found above, 
we conclude that:
\begin{equation*}
K^h_{tan}(x + t\vec n) =  \sqrt{e^h}(\mbox{Id} + th/h_0\Pi)^{-1}\Big[
2 t/h_0\, {\rm sym}(\nabla(A\vec n)-A\Pi)\Big](\mbox{Id} + th/h_0\Pi)^{-1} + o(\sqrt{e^h}).
\end{equation*}
Now, as $|\vec n_\e|=1$, the normal minor of $K^h$ is calculated by means of
(\ref{form2}) as:
\begin{equation*}
\vec n^T K^h(x+t\vec n)\vec n = |(\nabla_h y^h)\vec n|^2 - 1 = 
2th/h_0\e d^h \cdot \vec n_\e+ o(\sqrt{e^h})
= 2t/h_0 \sqrt{e^h} d^h \cdot \vec n +  o(\sqrt{e^h}).
\end{equation*}
The remaining coefficients of the symmetric matrix $K^h(x+ t\vec n)$ are, 
for $\tau\in T_x S$:
\begin{equation*}
\begin{split}
 \tau^T K^h(x+t\vec n)\vec n & = (\vec n_\e +th/h_0\e d^h)^T
\Big(\nabla u_\e+th/h_0\nabla \vec n_\e+
\frac{t^2}{2h_0^2}h^2\e\nabla d^h\Big)(\mbox{Id} + th/h_0\Pi)^{-1}\tau\\
& = t/h_0\sqrt{e^h}(d^h)^T\nabla u_\e(\mbox{Id} + th/h_0\Pi)^{-1}\tau + o(\sqrt{e^h}),
\end{split}
\end{equation*}
where we have used that $\vec n_\e^T\nabla\vec n_\e=0$.

{\bf 3.} From the previous computations we finally deduce that:
\begin{equation}\label{Kconv}
\lim_{h\to 0} \frac{1}{2\sqrt{e^h}} K^h = \frac{t}{h_0}K(x)_{tan}
+ \frac{t}{h_0}(\zeta\otimes \vec n + \vec n\otimes \zeta) 
\quad \mbox{ in } L^\infty(S^{h_0}),
\end{equation}
where the vector field $\zeta$ is defined in \eqref{def-zeta} and 
the symmetric matrix field $K_{tan}\in L^\infty(S, \mathbb{R}^{2\times 2})$ is given by:
\begin{equation}\label{Kdef}
K(x)_{tan} = {\rm sym} (\nabla(A\vec n) - A\Pi)_{tan}.
\end{equation}
Using \eqref{int-nv}, \eqref{Wdopp}, \eqref{Kconv} and the dominated
convergence theorem, we obtain:
\begin{equation*}
\begin{split}
\lim_{h\to 0} \frac{1}{e^h} E^h(u^h) & = 
\lim_{h\to 0} \frac{1}{e^h} \int_S \fint_{-h_0/2}^{h_0/2} W(\nabla_h y^h)
\det (\mbox{Id} + th/h_0\Pi)~\mbox{d}t\mbox{d}x\\
& = \frac{1}{2} \int_S \fint_{-h_0/2}^{h_0/2} \mathcal{Q}_3 \Big(\frac{t}{h_0}K(x)_{tan}
+ \frac{t}{h_0}(\zeta\otimes \vec n + \vec n\otimes \zeta)\Big)~\mbox{d}t\mbox{d}x \\
& = \frac{1}{2}\int_S \fint_{-h_0/2}^{h_0/2}\frac{t^2}{h_0^2}\mathcal{Q}_2
\big(x,{\rm sym} (\nabla(A\vec n) - A\Pi)_{tan}\big)~\mbox{d}t\mbox{d}x,
\end{split}
\end{equation*}
where the last equality is a consequence of \eqref{def-zeta} and
\eqref{Kdef}. Property (iii) now follows, upon integration with respect to $t$
in the last integral above.
\endproof

\medskip

\begin{remark}
Note that the recovery sequence formula (\ref{rec_seq}) follows the Cosserat
ansatz, which assumes that the fibers orthogonal to the mid-surface deform
linearly, only up to the term of order 1.  Naturally, the 
inhomogeneous stretch and twist (or warping) of these fibers 
is of order $\e$ and is due to the first order change in the second fundamental form of 
the surface, as follows from the formula (\ref{warp}) for the limit value of $d^h$. 
\end{remark}

\section{Convergence of minimizers in presence of body forces}
\label{sec_deadloads}

We recall that the major application of a $\Gamma$-convergence result, under
suitable compactness conditions, is that the minimizers of a given sequence of
functionals converge to the minimizers of their $\Gamma$-limit.
In view of Theorem \ref{thliminf-intro} and Theorem \ref{thlimsup-intro}, 
we have established the $\Gamma$-convergence of the rescaled elastic energy
$1/e^h E^h$ to the functional (\ref{I-intro}). 
In this section we shall see that under certain scaling 
regimes of forces, the scaled energies $1/e^h E^h$ of  minimizers of the total energies
$J^h$ in (\ref{total-intro}) are bounded, allowing then for an application 
of the compactness results in Theorem \ref{thliminf-intro}. 

\medskip

Let $f^h\in L^2(S^h,\mathbb{R}^3)$ be a sequence of forces acting on $S^h$,
of the form
$$f^h(x+t\vec n(x)) = h\sqrt{e^h}\det\left(\mbox{Id} + t\Pi(x)\right)^{-1} f^h(x),$$
where $f^h\in L^2(S,\mathbb{R}^3)$ have the properties:
\begin{equation}\label{fhass}
\int_S f^h = 0 \quad \mbox{and} \quad \lim_{h\to 0}f^h = f \mbox{ weakly in } L^2(S). 
\end{equation}
Let $m^h$  be the maximized action of force $f^h$ over all rotations
of $S$: 
$$m^h = \max_{Q\in SO(3)} \int_S f^h(x)\cdot Qx~\mbox{d}x,$$
and define: 
$$\mathcal{M}=\{\bar Q\in SO(3); ~ r(\bar Q)<+\infty\},$$ 
to be the effective domain of the following relaxation functional 
$r:SO(3)\longrightarrow [0,+\infty]$:
\begin{equation*}\label{rfunct}
r(Q) = \inf\left\{\liminf \frac{1}{e^h}\left(m^h - \int_S f^h(x)\cdot Q^h
    x\right);~~ Q^h\in SO(3), ~ Q^h\to Q\right\}.
\end{equation*}
In the particular case when $f^h=f$, one has $\mathcal{M} = \{\bar Q\in SO(3);
~ \int_S f\cdot \bar Q x= \max_{Q\in SO(3)}\int_S f\cdot Qx\}$.
As we shall see below, the set $\mathcal{M}$ identifies the candidates for
large rotations that the body chooses to perform in response to a force,
rather undergoing a further compression. 

The total energy functional on $S^h$ is given through:
$$J^h(u^h) = E^h(u^h)  - \frac{1}{h}\int_{S^h} f^h u^h + h\sqrt{e^h} m^h.$$
We then have the following result:

\begin{theorem}\label{thmaincinque}
Assume (\ref{scaling-intro}) and (\ref{fhass}). 
Let $S$ satisfy the ellipticity and regularity requirements stated in Theorem
\ref{thlimsup-intro}. Then:
\begin{itemize}
\item[(i)] For every small $h>0$ one has:
$$0\geq \inf\left\{\frac{1}{e^h} J^h(u^h); ~~ u^h\in W^{1,2}(S^h, \mathbb{R}^3)
\right\}\geq -C.$$
\item[(ii)]  If $u^h\in W^{1,2}(S, \mathbb{R}^3)$ is a minimizing sequence
of $\frac{1}{e^h} J^h$, that is:
$$ \lim_{h\to 0} \left(\frac{1}{e^h}J^h(u^h) - \inf\frac{1}{e^h}J^h
\right)=0, $$
then the conclusions of Theorem \ref{thliminf-intro} hold, and moreover any
accumulation point of $\{Q^h\}$ belongs to $\mathcal{M}$.
Further, any limit $(V, \bar Q)$ minimizes the functional: 
$$J(V,\bar Q) = I(V) - \int_S f\cdot \bar{Q} V + r(\bar Q),$$ 
over all $V\in \mathcal{V}$ and all $\bar Q\in \mathcal{M}$.
\end{itemize}
\end{theorem}
The proof follows exactly as in \cite{lemopa1, lemopa2}, hence we omit it. Notice that
when $f^h=f$, then the term $r(\bar Q)$ in the definition of the functional $J$ may
be dropped, as $r=0$ on $\mathcal{M}$.
In the general case, both $r$ and $\mathcal{M}$ depend on the asymptotic
behavior of the minimizers of the linear functions
$SO(3)\ni Q\mapsto \int_S f^h(x)\cdot Qx~\mbox{d}x.$     

For a further related discussion we refer to \cite{lemopa1}.

\end{document}